\documentclass[11pt]{article}
\usepackage{amssymb,epsfig,amsmath}

\usepackage{graphicx}
\usepackage{epstopdf}

\newtheorem{prop}{Proposition}
\newtheorem{defn}[prop]{Definition}
\newtheorem{lm}[prop]{Lemma}
\newtheorem{ex}[prop]{Example}
\newtheorem{thm}[prop]{Theorem}
\newtheorem{cor}[prop]{Corollary}

\newenvironment{pf}{\begin{trivlist}\item[] \textbf{Proof.}}
                     {\hspace*{\fill} $\square$\end{trivlist}}
                     
\newenvironment{pfl}[1]{\begin{trivlist}\item[] \textbf{Proof of Theorem #1.}}
                     {\hspace*{\fill} $\square$\end{trivlist}}

\newcommand{\beq}{\begin{eqnarray*}}
\newcommand{\eeq}{\end{eqnarray*}}

\newcommand{\bP}{\mathbb{P}}
\newcommand{\bN}{\mathbb{N}}
\newcommand{\bR}{\mathbb{R}}
\newcommand{\bE}{\mathbb{E}}
\newcommand{\bT}{\mathbb{T}}
\newcommand{\TT}{\widetilde{\mathcal{T}}}
\newcommand{\cC}{\mathcal{C}}
\newcommand{\cP}{\mathcal{P}}
\newcommand{\cS}{\mathcal{S}}
\newcommand{\cT}{\mathcal{T}}
\newcommand{\fs}{\mathbf{s}}
\newcommand{\ft} {\mathbf{t}}

\newcommand{\Ht}{\textrm{ht}}
\newcommand{\M}{\mathcal{M}_w}
\newcommand{\deq}{\mathrel{\mathop :}=}

\newcommand{\bs}{\mathbf{s}}

\newcommand*\samethanks[1][\value{footnote}]{\footnotemark[#1]}
\author{Foo\thanks{University of Podunk, Timbuktoo}
\and Bar\samethanks
\and Baz\thanks{Somewhere Else}
\and Bof\samethanks[1]
\and Hmm\samethanks}

\begin{document}

\title{
Regenerative tree growth: structural results and convergence\thanks{%
This research was supported by the National Science Foundation Awards 0405779 and 0806118, DR's work also in part by the National Science Foundation Graduate
Research Fellowship under Grant No. DGE 1106400 and in part by NSF DMS-1204840.}}
\author{
Jim Pitman%
\thanks{%
University of California at Berkeley; email pitman@stat.berkeley.edu, drizzolo@math.berkeley.edu} \and 
Douglas Rizzolo\samethanks \and%
Matthias Winkel\thanks{%
University of Oxford; email winkel@stats.ox.ac.uk}}

\maketitle


\begin{abstract} \noindent We introduce regenerative tree growth 
  processes as consistent families of random trees with $n$ labelled leaves, 
  $n\ge 1$, with a regenerative property at branch points. 
  This framework includes growth processes for exchangeably labelled Markov 
  branching trees, as well as non-exchangeable models such as the alpha-theta model, the alpha-gamma model and all restricted exchangeable models previously studied. Our main structural result is a representation of the growth rule 
  by a $\sigma$-finite dislocation measure $\kappa$ on the set of partitions of $\bN$ extending Bertoin's 
  notion of exchangeable dislocation measures from the setting of homogeneous fragmentations.
  We use this representation to establish necessary and sufficient conditions on the growth rule under 
  which we can apply results by Haas and Miermont for unlabelled and not necessarily consistent trees to establish self-similar 
  random trees and residual mass processes as scaling limits. While previous studies exploited some form of exchangeability,
  our scaling limit results here only require a regularity condition on the convergence of asymptotic frequencies under $\kappa$, in addition to a regular variation condition. 
    
\emph{AMS 2000 subject classifications: 60J80.\newline
}

\end{abstract}


\section{Introduction to regenerative tree growth processes}

For each $n\ge 1$, denote by $\bT_n$ the set of rooted leaf-labelled combinatorial trees with no degree-2 vertices and $n+1$ degree-1 vertices, one of which is called the \em root\em, the others \em leaves\em. We distinguish the leaves by labels $1,\ldots,n$. Vertices of degree 3 or higher are called \em branch points\em. Consider a family $\cT_n$, $n\ge 1$, of random trees in $\bT_n$, $n\ge 1$. For $n\ge 2$, we refer to the vertex adjacent to the root as the \em first branch point\em. It induces the \em first split\em, a random partition $\Pi_n=(\Pi_{n,1},\ldots,\Pi_{n,K_n})$ of the label set $[n]:=\{1,\ldots,n\}$ into the label sets of the \em subtrees above the branch point\em, the connected components of the tree with the first branch point removed. Here, we put the \em blocks \em $\Pi_{n,i}$ of $\Pi_n$ in the order of their least elements.
For illustration, we write\vspace{-0.2cm}
$$\bT_1=\left\{\parbox{0.27cm}{\includegraphics[scale=0.24]{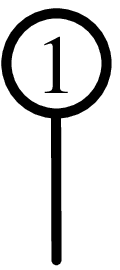}}\right\},\quad\bT_2=\left\{\parbox{0.6cm}{\includegraphics[scale=0.24]{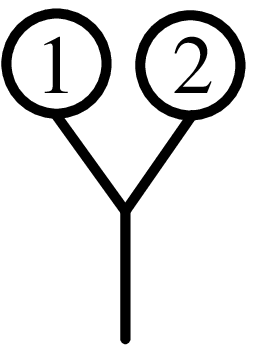}}\right\},\quad\bT_3=\left\{\parbox{0.9cm}{\includegraphics[height=1cm]{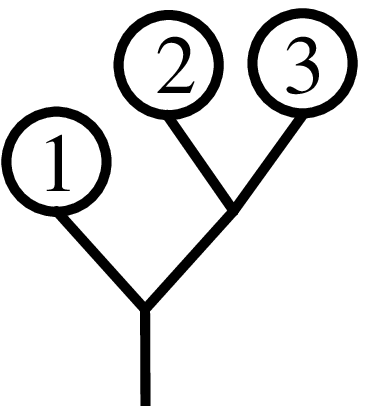}},\!\!\parbox{0.9cm}{\includegraphics[height=1cm]{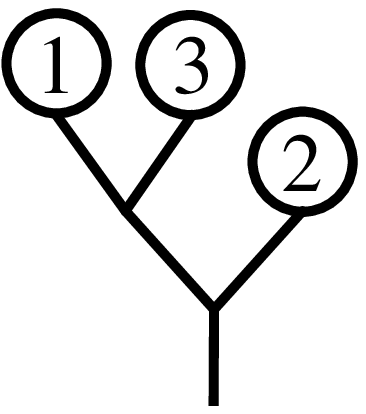}},\!\!\parbox{0.9cm}{\includegraphics[height=1cm]{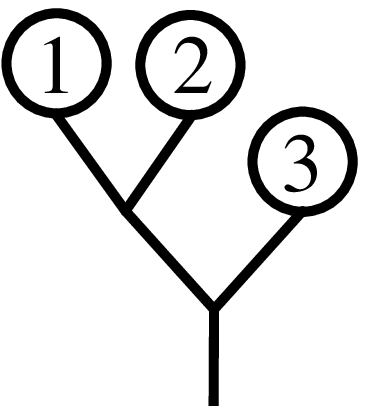}},\parbox{0.9cm}{\includegraphics[scale=0.24]{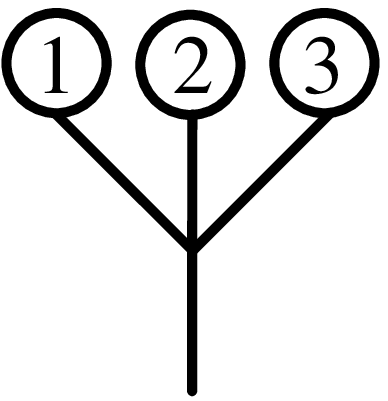}}\right\},\quad \mbox{etc.,}\vspace{-0.1cm}
$$
where we have ordered subtrees by their least labels to uniquely choose plane tree representatives.

We suppose that the family $(\cT_n,n\ge 1)$ is \em consistent \em in the sense that removal of leaf $n+1$ (and the resulting degree-2 vertex, if any) from $\cT_{n+1}$ yields $\cT_{n}$. Reversing this removal gives a tree growth step from $n$ to $n+1$. A consistent family $(\cT_n,n\ge 1)$ is called a \em tree growth process\em. 
For $B \subseteq [n]$, let $\bT_B$ be the set of trees with $\#B$ leaves 
  labelled by $B$, so that $\bT_{[n]}=\bT_n$. Let $\cT_{n,B} \in \bT_B$ \linebreak[1] be the \em reduced subtree \em of $\cT_n$ spanned by the root of $\cT_n$ and 
  leaves in $B$, and let $\TT_{n,B} \in \bT_{[\#B]}$ be the image of $\cT_{n,B}$ after \em relabelling of
  leaves \em by the increasing bijection from $B$ to $[\#B]$. 

\begin{defn}\label{df2} \rm We call a tree growth process $(\cT_n,n\ge 1)$ {\em regenerative} if for each $n \ge 2$, conditionally given that
  the first split of $\cT_n$ is $\Pi_n=(B_1,\ldots,B_k)$, the relabelled subtrees $\TT_{n,B_i}$, $1\le i\le k$, above the first branch point are independent copies of $\cT_{\#B_i}$. 
\end{defn}
While this property is well-known for many tree growth processes, the goals of the present paper are to provide general structural results and to study implications for continuum tree asymptotics in this general framework. In the terminology of \cite{CW}, the trees in a regenerative tree growth process as defined here are ``consistent labelled Markov branching trees''. The \em exchangeable \em case, where the distribution of $\cT_n$ is invariant under all permutations of labels, was initiated by Aldous \cite{Ald-93}, who posed the problem of providing a Kingman-type representation in this case. Bertoin's \cite{Ber-book} theory of homogeneous fragmentations solved that problem as explained in \cite{HMPW}. Then \cite{HMPW,HPW} studied tree growth processes associated with fragmentation processes. Natural non-exchangeable tree growth processes were described in terms of simple growth rules that admit regenerative descriptions based on the first split and its subtrees, see particularly \cite{CFW,For-05,PW09}, as reviewed in \linebreak[2] Examples \ref{ex1b} and \ref{ex1a} below.  We remark, however, that not all natural models fall into our current framework.  For example, if $\cT_n$ is uniform on $\bT_n$ then $(\cT_n,n\geq 1)$ is not a regenerative tree growth \linebreak[2] process because it is not consistent (see \cite{PR} for weak limits). Haulk and Pitman \cite{HaP-11} give de Finetti representations for exchangeable tree growth processes that are not necessarily regenerative.

\begin{figure}[t]
\begin{center}
\includegraphics[scale=0.35]{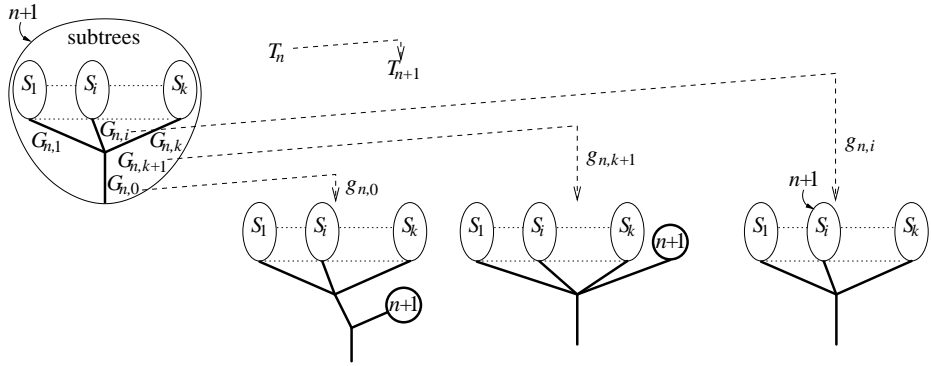}
\end{center}
\vspace{-0.8cm}
\caption{Illustration of the regenerative tree growth step}
\label{recgrowth}
\vspace{-0.2cm}
\end{figure}

An important consequence of Definition \ref{df2} is that all regenerative tree growth processes admit descriptions in terms of a growth rule (cf.\ Figure \ref{recgrowth}).

\begin{prop}\label{regfromg} In the tree growth step from $n$ to $n+1$ for $n\ge 2$, there are the following disjoint events, $G_{n,i}$ for $i=0,\ldots,K_n+1$, where $K_n\ge 2$ is the number of blocks of the first split of $\cT_n$:\vspace{-0.1cm}
  \begin{itemize}\item $G_{n,0}$: leaf $n+1$ is attached to a new branch point between the root and the first branch point of $\cT_n$;\vspace{-0.2cm}
    \item $G_{n,i}$, $1\le i\le K_n$: label $n+1$ is inserted into the $i$th block of the first split;\vspace{-0.2cm}
    \item $G_{n,K_n+1}$: leaf $n\!+\!1$ is attached to the first branch point, as singleton block of the first split.\vspace{-0.1cm}
  \end{itemize}
  A tree growth process $(\cT_n,n\ge 1)$ is regenerative if and only if $\bP(G_{n,0}\,|\,\cT_n)=\bP(G_{n,0})$ 
  does not depend on $\cT_n$ and $\bP(G_{n,i}\,|\,\cT_n)=\bP(G_{n,i}\,|\,\Pi_n)$, $1\le i\le K_n+1$, only 
  depends on the partition $\Pi_n$ of the first split. In the event $G_{n,i}$, $1\le i\le K_n$, label $n+1$  
  is inserted into the $i$th subtree of $\cT_n$ of size $\#\Pi_{n,i}$ following the same rule, up
  to relabelling by the increasing bijection from $\Pi_{n,i}\cup\{n+1\}$ to $[\#\Pi_{n,i}+1]$.  
\end{prop}
See Appendix \ref{ap_regfromg} for a proof of this proposition.  

We denote by $\cP_n$ the set of partitions $\pi=(B_1,\ldots,B_k)$ of $[n]$, with blocks $B_i$ ordered by their least elements. We use the notation $g_n(\pi,i)=\bP(G_{n,i}\,|\,\Pi_n=\pi)$, $0\le i\le k+1$, for $\pi\neq 1_{[n]}:=([n])$, $n\ge 2$, and write
$g_n(0)=g_n(\pi,0)$, since we require that this quantity is independent of $\pi\in\cP_n\setminus\{1_{[n]}\}$.

One of our main results is that regenerative tree growth rules are (almost) in one-to-one correspondence with $\sigma$-finite measures on $\cP$, the set of partitions of $\bN=\{1,2,3,\dots\}$. \pagebreak Before stating this, let us introduce the notation $\cP^\pi =\{\Gamma\in \cP : \Gamma^{[n]}=\pi\}$ where $\Gamma^{[n]}\in \cP_n$ is the partition whose blocks are the non-empty blocks of $(\Gamma_i\cap [n],i\geq 1)$.  Often we will abuse this notation and for a partition $\pi=(B_1,\dots, B_k)$ write $\cP^{B_1,\dots,B_k}$ instead of $\cP^\pi$.  The most common occurrence of this will be the use of $\cP^{[n]}$ instead of $\cP^{1_{[n]}}$.  We equip $\cP$ with the $\sigma$-algebra generated by 
$\{\cP^\pi,\pi\in\cP_n,n\ge 1\}$, which is also the Borel $\sigma$-algebra generated by the metric $d(\Gamma,\widehat{\Gamma})\!=\!\exp(-\inf\{n\!\ge\! 1\colon\Gamma^{[n]}\!\neq\!\widehat{\Gamma}^{[n]}\})$.  

\begin{thm} \label{thm_kappa1}
\begin{enumerate}
\item[{\em (i)}] Let $(g_n,n\geq 2)$ be a regenerative growth rule such that $g_j(0)<1$ for all $j\geq 2$.  Then there exists a unique $\sigma$-finite measure $\kappa$ on $\cP$ with $\kappa(\{1_{\bN}\})=0$ and $\kappa(\cP^{\{1\},\{2\}})=1$ such that \vspace{-0.2cm}
\begin{equation}\label{g_by_kappa}
g_n(0)=1-\frac{\lambda_n}{\lambda_{n+1}},\quad g_n(\pi,i)=\frac{\lambda_n}{\lambda_{n+1}}\frac{\kappa(\cP^{B_1,\ldots,B_{i-1},B_i\cup\{n+1\},B_{i+1},\ldots,B_k})}{\kappa(\cP^{B_1,\ldots,B_k})},\ i\in[k+1],\end{equation}
where $\lambda_n=\kappa(\cP\setminus\cP^{[n]})$ and $\pi = (B_1,\dots,B_k)$.  Moreover, in this case for $\pi\in \cP_n\setminus \{1_{[n]}\}$ we have the splitting rule $p_n(\pi):=\bP(\Pi_n=\pi) = \lambda_n^{-1}\kappa(\cP^\pi)$.

\item[{\em (ii)}] If $\kappa$ is any measure on $\cP$ such that $\kappa(\{1_\bN\})=0$, $0<\kappa(\cP^{\{1\},\{2\}})<\infty$, and $\kappa(\cP\setminus\cP^{[n]})<\infty$ for all $n\geq 2$, then $(g_n,n\geq 2)$ defined by {\rm \eqref{g_by_kappa}} is a regenerative growth rule such that $g_j(0)<1$ for all $j\geq 2$. In particular, there is a regenerative tree growth process associated with $\kappa$.  
\end{enumerate}
\end{thm}
We remark that part (ii) of this theorem shows how the relation between $(g_n,n\geq 2)$ and $\kappa$ fails to be one-to-one.  That is, if $\kappa$ produces a regenerative growth rule $(g_n,n\geq 2)$ by {\rm \eqref{g_by_kappa}}, then any constant multiple of $\kappa$ produces the same growth rule $(g_n,n\geq 2)$ by {\rm \eqref{g_by_kappa}}.  If $\kappa$ as in part (ii) and $(g_n,n\geq 2)$ are related by {\rm \eqref{g_by_kappa}}, we call $\kappa$ a \textit{dislocation measure} for $(g_n,n\geq 2)$.

Many of the asymptotic properties of a regenerative tree growth process can be obtained by analysing the asymptotic properties of the associated measure $\kappa$.  In fact, the two most important considerations turn out to be the growth rate of $\lambda_n$ and the regularity of the convergence of asymptotic frequencies under $\kappa$.  Let us expand on the second point.  For $\Gamma\in\cP$ and $n\ge 1$, consider the decreasing rearrangement 
$|\Gamma^{[n]}|^\downarrow\!=\!(|\Gamma^{[n]}|^\downarrow_i,i\!\ge\! 1)$ of relative frequencies
$(|\Gamma_i^{[n]}|,i\!\ge\! 1)$, where $|\Gamma_i^{[n]}|\!=\!\#\Gamma_i^{[n]}\!/n$.  If the limit as $n\rightarrow\infty$ of $|\Gamma^{[n]}|^\downarrow_i$ or of 
$|\Gamma^{[n]}_i|$ exists, this is denoted by $|\Gamma|^\downarrow_i$ and 
$|\Gamma_i|$, respectively, and we say that \em an asymptotic frequency exists \em for that part.  If $|\Gamma_i|$ exists for all $i$ we say $\Gamma$ has \textit{asymptotic frequencies} while if $|\Gamma|^\downarrow_i$ exists for all $i$ we say $\Gamma$ has \textit{asymptotic ranked frequencies}.  Moreover, if the asymptotic (ranked) frequencies exist and sum to $1$ $\kappa$-a.e., we say they are \textit{proper}.  If $\Gamma$ has asymptotic ranked frequencies then $|\Gamma|^\downarrow=(|\Gamma|^\downarrow_i, i\geq 1)$ naturally lives in the space\vspace{-0.5cm}
\begin{equation}\label{eq_sd}\cS^\downarrow = \left\{ (s_1,s_2,\dots)\colon s_1\geq s_2\geq \cdots \geq 0 \textrm{ and } \sum_{i\ge 1}s_i \leq 1\right\}.\vspace{-0.1cm}\end{equation}
We will equip $\cS^\downarrow$ with the topology of pointwise convergence (which is also the topology of $\ell_p$ convergence for any $p>1$).  We also introduce $\cS^\downarrow_1 = \{ \mathbf{s}\in\cS^\downarrow\colon \sum_{i\ge 1} s_i=1\}$.

We can then prove the following theorem, the background for which will be fully developed later.

\begin{thm}\label{thm_tree1}
Let $(\cT_n,n\geq 1)$ be a regenerative tree growth process associated with a dislocation measure $\kappa$.  Assume that $\kappa$-a.e.\ $\Gamma\in\cP$ has asymptotic ranked frequencies in $\cS^\downarrow_1\setminus\{(1,0,\ldots)\}$, define
      $\nu$ to be the push-forward of $\kappa$ under $\Gamma\mapsto|\Gamma|^\downarrow$ and suppose  $\int_{\cS^\downarrow}(1-s_1)\nu(d\fs)<\infty$ and $\lambda_n=\kappa(\cP\setminus\cP^{[n]})=n^\gamma\ell(n)$ for some slowly varying function $\ell$ and $\gamma>0$.  If  
\begin{equation}\label{eq_tree1} \lim_{n\to\infty}\int_{\cP}\left(|\Gamma^{[n]}|^\downarrow_1 -|\Gamma|^\downarrow_1\right)\kappa(d\Gamma) = 0,\end{equation}\vspace{-0.3cm}

\noindent then $\displaystyle\frac{\cT_n^\circ}{n^\gamma\ell(n)}\rightarrow\cT_{\gamma,\nu}$ in distribution, as $n\rightarrow\infty$, in the rooted Gromov-Hausdorff-Prokhorov (GHP) sense, where $\cT_{\gamma,\nu}$ is a self-similar fragmentation tree with characteristics $(\gamma,\nu)$ and $\cT^\circ_n$ is the tree obtained from $\cT_n$ by delabelling the leaves, considered as a metric measure space with the graph metric and the uniform probability measure on the leaves.
\end{thm}\pagebreak[2]
We remark that when considering $\cT_n^\circ$ purely as a tree we treat it as an element of the set $\bT_n^\circ$ of rooted unlabelled trees with $n$ leaves and no degree-2 vertices.

This theorem provides conditions for the existence of a scaling limit of $\cT^\circ_n$, where the label structure of $\cT_n$ has been forgotten.  However, the leaf labels are an integral part of the tree growth processes under consideration here, so it is natural to ask what happens to the labels.  
Ideally, one would like a notion of labelled continuum trees to serve as scaling limits of regenerative tree growth processes, just as there is a notion of ordered continuum trees that serve as scaling limits of ordered Galton-Watson trees \cite{Aldo93}.  
However, the appropriate notion is elusive, so we content ourselves with studying the leaf $\{1\}$ and the structure of the path from the root to this leaf.  We obtain several results relating the convergence of the residual mass process of the leaf $\{1\}$ to the existence of a scaling limit of the whole tree.  Here, by the residual mass process of $\{1\}$ we mean the Markov chain in 
$m\ge 0$ starting from $X^{(n)}_0=n$, decreasing to $X^{(n)}_1=\#\Pi_{n,1}$ and 
further according to successive splits of the block containing $\{1\}$ until $M_n=\inf\{m\ge 0\colon X^{(n)}_m=1\}$, when label 1 becomes a singleton.  We set $X^{(n)}_m=0$, $m>M_n$. The limiting processes are decreasing self-similar Markov processes in $[0,\infty)$, which Lamperti \cite{Lam-72} represented in terms of subordinators $\xi$, as
\begin{equation}\label{pssmp} X_t=\exp\left(-\xi_{\tau_\xi(t)}\right),\qquad\mbox{where}\quad\tau_{\xi}(t) = \inf\left\{ u\ge 0 : \int_0^u \exp\left( - \gamma\xi_s\right) ds> t\right\}.
\end{equation}

\begin{thm}\label{thm_mass1}
Let $(\cT_n,n\geq 1)$ be a regenerative tree growth process with dislocation measure $\kappa$ and $X^{(n)}$ the residual mass process of $\{1\}$ in $\cT_n$.  Assume that the first block $\Gamma_1$ of $\kappa$-a.e.\ $\Gamma\in\cP$ has an asymptotic frequency $|\Gamma_1|\in(0,1)$ and define $\Lambda$ as push-forward of $\kappa$ under $\Gamma\mapsto-\log(|\Gamma_1|)$. Suppose $\int_{(0,\infty)}(1-e^{-x})\Lambda(dx)<\infty$ and $\lambda_n=\kappa(\cP\setminus\cP^{[n]})=n^\gamma\ell(n)$ for some slowly varying function $\ell$ and $\gamma>0$.  If 
\begin{equation} \label{eq_mass1} \lim_{n\to\infty}\int_{\cP}\left(|\Gamma_1^{[n]}| -|\Gamma_1|\right)\kappa(d\Gamma)=0
\end{equation} 
then $X^{(n)}_{\lfloor\lambda_n t\rfloor}/n\rightarrow X_t$ in distribution, as $n\!\rightarrow\!\infty$, in the Skorohod sense as functions of $t\ge 0$, where $X$ is a self-similar Markov process and  $\bE(e^{-s\xi_r})=\exp(-r\int_{(0,\infty)}(1\!-\!e^{-sy})\Lambda(dy))$ in Lamperti's representation {\rm(\ref{pssmp})}.  
Moreover, letting $A_n$ be the absorption time of $X^{(n)}$ at $0$, the above convergence in distribution holds jointly with the convergence of $A_n/\lambda_n$ to the absorption time of $X$ at $0$.  

  If, in addition, $\kappa$-a.e.\ $\Gamma\in\cP$ has asymptotic ranked frequencies then $\displaystyle\frac{\cT_n^\circ}{n^\gamma\ell(n)}\rightarrow\cT_{\gamma,\nu}$ in distribution, as $n\rightarrow\infty$, in the rooted GHP sense, as in Theorem \ref{thm_tree1}.
\end{thm}
Assuming that $\kappa$-a.e.\ $\Gamma\in\cP$ has asymptotic ranked frequencies, the remaining conditions of Theorem \ref{thm_mass1} are stronger than those of Theorem \ref{thm_tree1}.  In particular, note that
\[\int_{\cS^\downarrow}(1-s_1)\nu(d\fs) \leq \int_{(0,\infty)}(1-e^{-x})\Lambda(dx)
\]
and that \eqref{eq_mass1} implies \eqref{eq_tree1} (see the proof of Theorem \ref{propmc}).
In Example \ref{ex15} we construct a regenerative tree growth process where the conditions of Theorem \ref{thm_tree1} are satisfied, but the conditions of Theorem \ref{thm_mass1} are not.  We note again that leaf $\{1\}$ is generally not typical (i.e.\  uniformly random) and a heuristic interpretation of the last part of Theorem \ref{thm_mass1} is that the natural conditions that imply the convergence of the residual mass process of leaf $\{1\}$ are strong enough to imply that the residual mass process of a typical leaf converges as well.

In the other direction, there is a natural strengthening of the hypotheses of Theorem \ref{thm_tree1} that implies that the conclusions of Theorem \ref{thm_mass1} are satisfied.

\begin{cor}\label{cor p-con}
If, in addition to the hypotheses of Theorem \ref{thm_tree1}, including {\rm(\ref{eq_tree1})}, we assume that\vspace{-0.05cm}
\[\kappa(|\Gamma_1|\neq|\Gamma|^\downarrow_1)<\infty \quad \textrm{and} \quad \lim_{n\to\infty} \int_{\{|\Gamma_1|=|\Gamma|^\downarrow_1\}}\left(|\Gamma^{[n]}|^\downarrow_1 -|\Gamma_1^{[n]}| \right)\kappa(d\Gamma)=0\vspace{-0.15cm}\]
then, with the notation of Theorem \ref{thm_mass1}, $X^{(n)}_{\lfloor\lambda_n t\rfloor}/n\rightarrow X_t$ in distribution, as $n\!\rightarrow\!\infty$, in the Skorohod sense as functions of $t\ge 0$ and this convergence in distribution happens jointly with the convergence of $\lambda_n^{-1}A_n$ to the absorption time of $X$ at $0$.
\end{cor}
When $\kappa$-a.e.\ $\Gamma\in\cP$ has asymptotic ranked frequencies, Theorem \ref{thm_mass1}, combined with previous results about the residual mass process of a typical leaf, provides a description of how leaf $\{1\}$ differs from a typical leaf.  To see this, let $(U^{(n)}_k , k\geq 0  )$ be the residual mass process of a leaf picked uniformly at random from $\cT_n$ for $n\geq 1$.  Under the assumptions of Theorem \ref{thm_tree1}, Lemma 28 in \cite{HM10} implies that $U^{(n)}_{\lfloor\lambda_n t\rfloor}/n\rightarrow U_t=\exp(-\zeta_{\tau_\zeta(t)})$ where $(\zeta_t,t\geq 0)$ is a subordinator with 
\[\bE(e^{-s\zeta_r})= \exp\left(- r \int_{\cP} \left(1-\sum_{i\geq 1} (|\Gamma|^{\downarrow}_i)^{s+1} \right)\kappa(d\Gamma)\right).\]
It is easy to check that this agrees with the expression for $\bE (e^{-s\xi_r})$ in Theorem \ref{thm_mass1} when $\kappa$ is exchangeable (see Example \ref{ex2}), but these two expressions may differ in general.  Thus we have identified the scaling limit of the residual mass process of $\{1\}$ and scaling limit of the residual mass process of a uniform leaf in terms of subordinators whose Laplace exponents we know explicitly in terms of $\kappa$.  This provides insight into the difference between what the tree looks like from $\{1\}$'s perspective versus that of a typical leaf.

The remainder of this paper is organized as follows.  In Section \ref{sec_structure} we give a detailed analysis of the laws of $\cT_n$.  The proof of Theorem \ref{thm_kappa1} can be found here, along with a number of other structural results.  Section \ref{sec_examples} is devoted to examples.  One of the nice aspects of the theory presented in this paper is that it gives a coherent framework for many particular models that have been studied previously in the literature.  As a result, we are able to give simplified proofs of a number of previously known results about these models.  Moreover, our framework makes it easy to specify regenerative growth processes with desired asymptotic properties and this allows us to construct examples illustrating what can go wrong if some of the hypotheses of our theorems are left out.  Section \ref{sec_background} provides the necessary background to understand the precise meaning of our statements about scaling limits.  We define the limit objects $\cT_{\gamma,\nu}$, the GHP topology, and provide the main results from the literature on which our present theorems are built.  In Section \ref{subsecHMtrees} and \ref{secHMmcs} respectively, we provide the proofs of Theorems \ref{thm_tree1} and \ref{thm_mass1} based on general convergence criteria by Haas and Miermont \cite{HM09, HM10} for (not necessarily consistent) Markov branching models and non-increasing Markov chains.  Actually, our results are stronger but a bit more technical than Theorems \ref{thm_tree1} and \ref{thm_mass1}, so we will prove results that have these theorems as obvious consequences. Section \ref{further} gives some pointers at further problems and related work.

\section{Laws of regenerative growth processes}\label{sec_structure}
\subsection{Explicit formulas in terms of the growth rule}
The regenerative nature of the growth processes conditioned on the partition at the first split shows that much of the analysis of these processes can be reduced to analyzing the laws of the partition at the first split of $\cT_n$, i.e.\ the \textit{splitting rule} $p_n$, $n\ge 2$. We first find the splitting rule in terms of $(g_n,n\geq 2)$ and then obtain a formula for the law of $\cT_n$.  From the growth rule, we have for all $\pi=(B_1,\ldots,B_k)\in\cP_n\setminus\{1_{[n]}\}$,
\begin{equation}\label{pg}\begin{array}{c}p_2(\{1\},\{2\})=1,\qquad p_{n+1}([n],\{n+1\})=g_n(0),\quad n\ge 2,\\[0.2cm]
p_{n+1}(B_1,\ldots,B_{i-1},B_i\cup\{n+1\},B_{i+1},\ldots,B_k)=p_n(B_1,\ldots,B_k)g_n(\pi,i),\qquad i\in[k+1].\end{array}
\end{equation}
Using the natural convention $g_1(0)=1$, the solution to these equations can be written as\vspace{-0.1cm}
\begin{equation}\label{splrules}p_n(\pi)=p_n(B_1,\ldots,B_k)=g_{\min B_2-1}(0)\prod_{j=\min B_2}^{n-1} g_j(\pi^{[j]},I_j),\quad\mbox{where $I_j=i$ if $j+1\in B_i$,}\vspace{-0.1cm}\pagebreak
\end{equation}
and
$\pi^{[j]}$ is the vector of non-empty $B_i\cap[j]$. The RHS of this formula
is the probability of successively creating a new first branch point when $\min B_2$ is added and inserting all higher labels such that the resulting partition at the first split is $\pi$. By 
the regenerative property of $\cT_n$, we can write tree probabilities as a product over branch points; for a tree $T\in\bT_n$, we identify each vertex with the set $B$ of labels in the subtrees above this vertex, write $\pi(B)$ for the partition of the split at $B$, and $\widetilde{\pi}(B)$ for the partition
of $[\#B]$ obtained when relabelling $\pi(B)$ by the increasing bijection from $B$ to $[\#B]$:\vspace{-0.2cm}
\begin{equation}\label{treeprob}\bP(\cT_n=T)=\!\!\!\prod_{B\in T:\#B\ge 2}\!\!p_{\#B}(\widetilde{\pi}(B))=\!\!\!\prod_{B\in T:\#B\ge 2}\!\!\left(g_{\min\widetilde{\pi}(B)_2-1}(0)\!\!\!\prod_{j=\min\widetilde{\pi}(B)_2}^{\#B-1}\!\! g_j(\widetilde{\pi}(B)^{[j]},I_j(B))\right)\!\!,\end{equation}
where $I_j(B)=i$ if $j+1\in\widetilde{\pi}(B)_i$, and where $\widetilde{\pi}(B)=(\widetilde{\pi}(B)_1,\ldots,\widetilde{\pi}(B)_{k(B)})$. 

The residual mass process of the leaf $\{1\}$ can be described in terms of $p_n$.  Recall that the residual mass process is a Markov chain in 
$m\ge 1$ starting from $X^{(n)}_0=n$, decreasing to $X^{(n)}_1=\#\Pi_{n,1}$ and 
further according to successive splits until $M_n=\inf\{m\ge 0\colon X^{(n)}_m=1\}$, when label 1 becomes a singleton. We represent this Markov chain as a 
composition of $n$
$$\cC_n=\left(C^{(n)}_0,\ldots,C^{(n)}_{M_n}\right)=\left(X^{(n)}_0-X^{(n)}_1,X^{(n)}_1-X^{(n)}_2,\ldots,X^{(n)}_{M_n-1}-X^{(n)}_{M_n},X^{(n)}_{M_n}\right).$$
\begin{prop}\label{resmp} In a regenerative tree growth process, the family $(\cC_n,n\ge 1)$ of compositions is regenerative in the sense
  that conditionally given $C^{(n)}_0\!=j$, the composition $(C^{(n)}_1\!,\ldots,C^{(n)}_{M_n})$ of $n-j$ has the same distribution as $\cC_{n-j}$. The entries of the transition probability matrix are 
  $$\bP(C^{(n)}_0=n-j)=\bP(X^{(n)}_1=j)=\sum_{\pi=(B_1,\ldots,B_k)\in\cP_n\colon\#B_1=j}p_n(\pi),\qquad 1\le j \le n-1.$$
\end{prop}
This is a straightforward consequence of Definition \ref{df2}. We stress that we have consistency in the sense that $\cC_n$ can be obtained from $\cC_{n+1}$ by reducing one part of $\cC_{n+1}$ by 1 (the one
corresponding to label $n\!+\!1$ in $\cT_{n+1}$), but $(\cC_n,n\!\ge\! 1)$ is not
\em sampling consistent \em in the sense of \cite{gp} as this part is not a size-biased pick from $\cC_{n+1}$, in general. In special cases, versions of
this proposition are in the literature; in the exchangeable (sampling consistent) case, it is implicit in Bertoin's \cite{Ber-book} study of tagged particles and explicit in \cite{HPW}.

\subsection{The dislocation measure}
Recall the notation $\cP_n$ for the set of partitions $\pi\!=\!(B_1,\ldots,B_k)$ of 
$[n]\!=\!\{1,\ldots,n\}$ with blocks indexed in increasing order of their least elements and the notation $\cP$ for the set of all partitions $\Gamma=(\Gamma_i,i\ge 1)$ of $\bN$, with blocks ordered by their least element and $\Gamma_i=\varnothing$ if there are fewer than $i$ blocks.  Theorem \ref{thm_kappa1} relates growth rules $g_n$ and splitting rules $p_n$ on $\cP_n\setminus\{1_{[n]}\}$ to dislocation measures $\kappa$ on $\cP$.

Taking our cues from the exchangeable case, cf.\ \cite{Ber-book}, one thing we want from our dislocation measures is to be able to use them to embed regenerative tree growth processes in continuous time, making the trees the genealogical trees of continuous-time fragmentation processes.  This $\kappa$ is to provide rates $$\lambda_n=\kappa(\{\Gamma\in\cP\colon\Gamma^{[n]}\neq 1_{[n]}\})=\kappa(\cP\setminus\cP^{[n]})$$
for the first split of $[n]$, $n\ge 2$, which allow us to consistently embed the evolution of blocks in $\cT_n$, $n\ge 1$, into continuous time (see Theorem \ref{homf}).  Observe that the rate $\lambda_n$ of the first split of $[n]$ can then be thinned by the event that this split also
splits $[n-1]$, an event with probability $1-g_{n-1}(0)$, where $(g_n,n\ge 2)$ is the growth rule of the regenerative tree growth process, so that we need
\begin{equation}\label{lambdas}\lambda_n(1-g_{n-1}(0))=\lambda_{n-1},\quad n\ge 3,\quad\mbox{and hence }\lambda_n=\lambda_2\prod_{j=2}^{n-1}\frac{1}{1-g_j(0)},\quad\mbox{if $g_j(0)\neq 1$, $j\ge 2$.}\pagebreak\end{equation}
Note that $g_j(0)=1$ for any $j\ge 2$ means that all insertions in a subtree with $j$ leaves are made below the first split; if scaling limits of $(\cT_n,n\ge 1)$ exist at all, such subtrees with $j$ leaves will collapse in 
the scaling as $n\rightarrow\infty$. We will exclude such behaviour in the sequel and make the
\bigskip

\centerline{\noindent\bf Assumption (A) \rm $\ g_j(0)<1$ for all $j\ge 2$.} 

\begin{prop}\label{prop6} Consider a regenerative tree growth rule $(g_n,n\ge 2)$ 
  satisfying Assumption {\bf(A)} with splitting rule $(p_n,n\ge 2)$ given by {\rm(\ref{splrules})}, and let $\lambda_2>0$ be arbitrary. With $\lambda_n$, $n\ge 3$, defined by {\rm(\ref{lambdas})}, define\vspace{-0.3cm}
\begin{equation}\kappa(\cP^\pi)=\lambda_np_n(\pi),\qquad\pi\in\cP_n\setminus\{1_{[n]}\},n\ge 2;\qquad \kappa(\{1_\bN\})=0.\label{kappa}\end{equation}
Then $\kappa$ extends uniquely to a measure on $\cP$.\pagebreak[2] 
\end{prop}

\begin{pf} This is essentially the same as the analogous result for exchangeable fragmentations, cf.\ Bertoin's argument \cite[Proposition 3.2]{Ber-book}.  We have defined $\kappa$ on $\{\cP^\pi, \pi \in \bigcup_n(\cP_n\setminus \{1_{[n]}\})\}\cup \{1_{\bN}\}$ and it clearly extends to a countably additive measure on the ring generated by these sets.   Carath\'eodory's Extension Theorem provides the unique extension to the $\sigma$-ring these sets generate in $\cP$.  It is then straightforward to check that this $\sigma$-ring is a $\sigma$-algebra and, in fact, is the Borel $\sigma$-algebra on $\cP$.
\end{pf} 

Note that we can condition $\kappa$ on splitting $[n]$ and write $p_n$ 
as \begin{equation}p_n(\pi)=\kappa(\cP^\pi)/\kappa(\cP\setminus\cP^{[n]}),\qquad\pi\in\cP_n\setminus\{1_{[n]}\},\ n\ge 2.\label{pkappa}\end{equation}

\begin{pfl}{\ref{thm_kappa1}}
Theorem \ref{thm_kappa1} is a direct consequence of \eqref{pg} and Proposition \ref{prop6}.
\end{pfl}
Let us discuss how Bertoin's \cite{Ber-book} notion of a 
$\cP$-valued homogeneous fragmentation process finds a natural extension where his exchangeable dislocation 
measure is replaced by a dislocation measure in the sense defined above. 

\begin{defn}\rm A $\cP$-valued process $\Pi=(\Pi(t),t\ge 0)$ is called \em refining \em if for all $s<t$ and all blocks $\Pi_j(t)$ of $\Pi(t)$, there is a block $\Pi_i(s)$ of $\Pi(s)$ that contains $\Pi_j(t)$. For
a refining process $\Pi$, we define \em genealogical trees \em $\cT_n\in\bT_n$, $n\ge 1$, using the representation above (\ref{treeprob}): $\cT_n$ has as branch points and leaves all blocks $\Pi_i^{[n]}(t)$, $i\ge 1$, $t\ge 0$, visited by the restriction $\Pi^{[n]}$ of $\Pi$ to $[n]$.
\end{defn} 
Every regenerative tree growth process can be represented by a nice refining $\cP$-valued process:

\begin{thm}\label{homf} For each dislocation measure $\kappa$ as defined after Theorem \ref{thm_kappa1}, there exists a $\cP$-valued Feller process $\Pi=(\Pi(t),t\ge 0)$ such that the genealogical trees $\cT_n$ of the restrictions $\Pi^{[n]}$ of $\Pi$ to $[n]$, $n\ge 1$, form a regenerative tree growth process associated with dislocation measure $\kappa$.
\end{thm}
\begin{pf} We will use $\kappa$ in a Poissonian construction based on
independent $\cP$-valued Poisson point processes $(\Xi^{(i)}(t),t\ge 0)$, $i\ge 1$, with intensity measure $\kappa$. 
Roughly, we construct $\Pi$ with $\Pi(0)=1_\bN$ such that for all $i$ and $t$ 
the partition $\Xi^{(i)}(t)$ fragments the $i$th block $\Pi_i(t)$ of $\Pi(t)$ into the image 
$\widetilde{\Xi}^{(i)}(t)$ of $\Xi^{(i)}(t)$ under the increasing bijection from $\bN$, or $[\#\Pi_i(t)]$, to $\Pi_i(t)$.

More precisely, we build consistent $\cP_n$-valued continuous-time Markov chains $(\Pi^{[n]}(t),t\ge 0)$, $n\ge 1$, with jump times $S^{[n]}(k)\ge 0$ and states $M^{[n]}(k)=(M_1^{[n]}(k),\ldots,M_{K^{[n]}(k)}^{[n]}(k))\in\cP_n$: we set $\Pi^{[n]}(t)=M^{[n]}(k)$, $S^{[n]}(k)\le t<S^{[n]}(k+1)$, $k\ge 0$, where $S^{[n]}(0)=0$, $M^{[n]}(0)=1_{[n]}$, 
$$ S^{[n]}(k+1)=\inf\left\{t>S^{[n]}(k)\colon [\#M_i^{[n]}(k)]\not\in\Xi_1^{(i)}(t)\mbox{ for some $1\le i\le K^{[n]}(k)$}\right\}$$
and, if $S^{[n]}(k+1)<\infty$, let $M^{[n]}(k+1)$ be the partition obtained from $M^{[n]}(k)$ by replacing the $i$th block by the blocks of $\widetilde{\Xi}^{(i)}(S^{[n]}(k+1))$, the image of $\Xi^{(i)}(S^{[n]}(k+1))\cap[\#M_i^{[n]}(k)]$ under the increasing bijection from $[\#M_i^{[n]}(k)]$ to $M_i^{[n]}(k)$. Note that $S^{[n]}(k+1)=\infty$ if and only if $M^{[n]}(k)=0_{[n]}:=(\{1\},\ldots,\{n\})$,
as we require $\lambda_2=\kappa(\cP^{\{1\},\{2\}})>0$ for all dislocation measures.

Since $\Pi$ is uniquely determined by $(\Pi^{[n]},n\ge 1)$, standard properties
of Poisson point processes, and of the space $\cP$ complete the proof.
\end{pf}
By Theorem \ref{thm_kappa1}, a growth rule $(g_n,n\ge 2)$ determines a measure $\kappa$ only
up to a multiplicative factor $\lambda_2>0$. This is reflected in the 
fragmentation processes $\Pi$ of Theorem \ref{homf} in the fact 
that the genealogical trees $\cT_n$, $n\ge 1$, are unaffected by (linear) time changes of $\Pi$. \pagebreak[1]

From the consistency of $(\cT_n,n\ge1)$, it is clear that there is a unique 
branch point of $\cT_n$, where 1 and 2 are separated into different blocks. 
Moreover, as $n$ varies, the partitions at this branch point define a partition of 
some random subset of $\bN$, whose distribution when relabelled by the 
increasing bijection is described by the splitting rule conditioned on partitions 
that 
restrict to $(\{1\},\{2\})$, hence by 
$\kappa(\,\cdot\,|\,\cP^{\{1\},\{2\}})$. In the Poissonian construction, this 
partition after relabelling is $\Xi^{(1)}(S^{[2]}(1))$. 
More generally, while there may be Poisson points that do not induce branch 
points of $\cT_n$, $n\ge 1$, e.g.\ when $\kappa$ is finite or when $\kappa$ can produce blocks
of finite size, those points $\Xi^{(i)}(S^{[n]}(k))$ used in the Poissonian 
construction describe the partition at a branch point of $\cT_m$ for all $m\ge n$.
The partition at every branch point, separating labels $j$ and $\ell$ say, has  
a distribution that is absolutely continuous with respect to $\kappa$. 

The Poissonian construction formulated here differs from Bertoin's 
\cite[Section 3.1.3]{Ber-book} in the relabelling by increasing bijections: Bertoin uses
$\Pi_i(t)\cap\Xi^{(i)}(t)$ instead of $\widetilde{\Xi}^{(i)}(t)$. In the 
exchangeable case this yields the same processes, in distribution. A notable consequence is 
that under assumptions that ensure that there are always infinitely many blocks
and that they are all infinite, we can recover the $(\Xi^{(i)},i\ge 1)$  
from $\Pi$ in our setting. It is now possible to explore some more of Bertoin's fragmentation theory
\cite[Chapter 3]{Ber-book} in our extended generality, notably erosion
effects, stopping lines and extended branching properties, and, under conditions that ensure the existence of asymptotic frequencies, also self-similar
partition-valued fragmentation processes. More generally, it would be interesting to characterise Markov processes (with a suitable branching property) whose genealogical trees are regenerative. \vspace{-0.2cm}

\subsection{Unlabelled Markov branching trees}\label{subsec M-branch}
Our scaling limit results take advantage of recent progress on scaling limits of rooted unlabelled Markov branching trees, which we now introduce.  For $n\geq 1$, let $\bT^\circ_n$ be the image of $\bT_n$ under the map that delabels the leaves of a tree.  Define\vspace{-0.1cm}
\[\cP_n^\circ = \Bigg\{ (n_1,\dots, n_p)\in \bigcup_{k\geq 1} \bN^k : n\ge n_1\geq n_2\geq \cdots\geq n_p  \ \textrm{ and } \ \sum_{i=1}^pn_i=n, p\ge 1\Bigg\}.\vspace{-0.1cm}\]
Let $(p_n^\circ,n\geq 2)$ be a sequence such that for each $n$, $p_n^\circ$ is a probability function on $\cP_n^\circ\setminus\{(n)\}$.  A sequence $(\cT^\circ_n,n\geq 1)$ of random variables such that $\cT^\circ_n\in \bT^\circ_n$ is called a \textit{Markov branching model} based on $(p_n^\circ,n\geq 2)$ if for each $n\geq 2$, the law of $\cT^\circ_n$ is the same as the law of the tree $\widehat\cT$ constructed as follows: choose $(N_1,\dots,N_p)$ according to $p_n^\circ$; conditionally given that $(N_1,\dots,N_p)=(n_1,\dots, n_p)$, let $(\widehat{\cT_1},\dots, \widehat{\cT_p})$ be a vector of independent trees such that $\widehat{ \cT_i}$ is distributed as $\cT^\circ_{n_i}$, $1\le i\le p$; the tree $\widehat{\cT}$ is then obtained by identifying the roots of $\widehat{\cT_1},\dots, \widehat{\cT_p}$ as a single vertex and attaching a new root to this vertex. 

The following proposition is an immediate consequence of these definitions.

\begin{prop}\label{prop_mb}
If $(\cT_n,n\geq 1)$ is a regenerative tree growth process with associated dislocation measure $\kappa$ and $(\cT^\circ_n,n\geq 1)$ is the sequence of trees such that $\cT^\circ_n$ is obtained from $\cT_n$ by delabelling the leaves, then $(\cT^\circ_n,n\geq 1)$ is a Markov branching model based on the functions\vspace{-0.1cm}
\[ p_n^\circ(n_1,\ldots,n_k)=\sum_{\pi\in\cP_n\colon(\#\pi)^\downarrow=(n_1,\ldots,n_k)}p_n(\pi)=\sum_{\pi\in\cP_n\colon(\#\pi)^\downarrow=(n_1,\ldots,n_k)}\frac{\kappa(\cP^\pi)}{\lambda_n},\quad \lambda_n=\kappa(\cP\setminus\cP^{[n]}),\]
where we write $(\#\pi)^\downarrow$ for the decreasing rearrangement of the block sizes of $\pi$.
\end{prop}\vspace{-0.2cm}

\section{Examples}\label{sec_examples}\vspace{-0.1cm}

An important motivation for our results is that they allow a unified treatment of previously studied models. In this section we discuss these models, recall or construct their dislocation \pagebreak measures and demonstrate how our Theorems \ref{thm_tree1} and \ref{thm_mass1} apply. We also develop some further examples that explore the conditions (\ref{eq_tree1}) and (\ref{eq_mass1}) that appear in Theorems \ref{thm_tree1} and \ref{thm_mass1}. Before proceeding with the examples, we introduce paintbox partitions, which are a recurring theme in the construction of dislocation measures.  For $\fs\in\cS^\downarrow$, where $\cS^\downarrow$ is defined in \eqref{eq_sd}, we define \em Kingman's paintbox \em $\kappa_\fs$  as the distribution of the random partition $\Pi$ of $\bN$ where $i,j\in\bN$ are in the same block if $i=j$ or $R_i=R_j\ge 1$, where the $R_i$, $i\in\bN$, are independent random variables with $\bP(R_i=k)=s_k$, $k\ge 0$, and where $s_0=1-\sum_{i\ge 1}s_i$. Note that the Strong Law of Large Numbers implies that $\kappa_\fs$-a.e.\ $\Gamma\in\cP$ has asymptotic ranked frequencies $|\Gamma|^\downarrow=\fs$.

\begin{ex}[Exchangeable models \cite{Ber-book,HMPW}]\label{ex2}\rm  Bertoin classified all
  \em exchangeable \em dislocation measures, i.e.\ measures that are invariant under the action of permutations of $\bN$ on $\cP$, giving an integral representation\vspace{-0.3cm}
  $$\kappa=\sum_{j\ge 1}c\delta_{\varepsilon^{(j)}}+\int_{\cS^\downarrow}\kappa_\fs(\,\cdot\,)\nu(d\fs),\vspace{-0.1cm}$$
  where $c\ge 0$, $\varepsilon^{(j)}$ is the partition with blocks $\{j\}$ and $\bN\!\setminus\!\{j\}$, and $\nu$ is a measure on $\cS^\downarrow$ with
  \begin{equation}\label{condnu1}\nu(\{(1,0,0,\ldots)\})=0\quad\mbox{and}\quad\int_{\cS^\downarrow}(1-s_1)\nu(d\fs)<\infty.
  \end{equation}
  Then $\nu$ is the push-forward of $\kappa$ under  
  $\Gamma\mapsto|\Gamma|^\downarrow$, restricted to $\cS^\downarrow\setminus\{(1,0,\ldots)\}$.

The splitting rules $(p_n,n\ge 2)$ associated with Bertoin's exchangeable dislocation measures $\kappa$ via \eqref{kappa} give rise to the consistent exchangeably 
  labelled Markov branching trees of \cite{HMPW}. For exchangeable  
  $\kappa=\int_{\cS^\downarrow_1}\kappa_\fs(\cdot)\nu(d\fs)$, it was
  demonstrated in \cite{HM10} that if $\lambda_n=\kappa(\cP\setminus\cP^{[n]})=n^\gamma\ell(n)$ for some $\gamma>0$ and some slowly varying function $\ell$, then \eqref{eq_tree1} also holds.  It is also easy to verify the condition of Theorem \ref{thm_mass1} in this case.  
\hfill $\square$ \end{ex}
Example \ref{ex2} includes an important subclass of models previously studied in their own right, whose dislocation measures are of Poisson-Dirichlet type. For this example, we can also calculate the growth rule explicitly up to a sequence of normalisation constants:

\begin{ex}[Poisson-Dirichlet model \cite{HPW,MPW}]\label{ex3}\rm According to \cite{MPW}, the only consistent exchangeable model with splitting rule of the Gibbs form\vspace{-0.1cm} 
$$p_n(\pi)=\frac{a_k}{c_n}\prod_{i=1}^kw_{\#B_i},\quad 
 \pi=(B_1,\ldots,B_k)\in\cP_n\setminus \{1_{[n]}\},\qquad\mbox{for some $w_j\ge 0$, $a_k\ge 0$, $c_n>0$,}\vspace{-0.1cm}$$
is given by a two-parameter family. Most relevant for us are $0<\alpha<1$ and $\theta\ge-2\alpha$ with $w_j=\Gamma(j-\alpha)/\Gamma(1-\alpha)$, $j\ge 1$, and  $a_k=\alpha^{k-2}\Gamma(k+\theta/\alpha)/\Gamma(2+\theta/\alpha)$, $k\ge 2$, with normalisation constants  $c_n=c_{\alpha,\theta}(n)$ satisfying $c_{\alpha,\theta}(2)=1$ and  $c_{\alpha,\theta}(n+1)=(n+\theta)c_{\alpha,\theta}(n)+\Gamma(n-\alpha)/\Gamma(1-\alpha)$, $n\ge 2$. Case $\alpha=0$ is a limiting case. 
 These yield growth rules for $\pi=(B_1,\ldots,B_k)$ of the form
  \beq g_n(0)\!\!&=&\!\!p_{n+1}([n],\{n+1\})=\frac{\Gamma(n-\alpha)}{\Gamma(1-\alpha)c_{\alpha,\theta}(n+1)},\\[0.1cm]
	   g_n(\pi,i)\!\!&=&\!\!\frac{p_{n+1}(B_1,\ldots,B_{i-1},B_i\cup\{n+1\},B_{i+1},\ldots,B_k)}{p_n(B_1,\ldots,B_k)}=\frac{(\#B_i-\alpha)c_{\alpha,\theta}(n)}{c_{\alpha,\theta}(n+1)},\quad i\in[k],\\[0.1cm]
       g_n(\pi,k+1)\!\!&=&\!\!\frac{(k\alpha+\theta)c_{\alpha,\theta}(n)}{c_{\alpha,\theta}(n+1)}.
  \eeq
Dislocation measures $\kappa_{\alpha,\theta}^{\rm PD^*}$ are exchangeable with, for $\theta>-2\alpha$, $\alpha\in(0,1)$, coefficient $c=0$ and\vspace{-0.1cm}
     $$\nu(d\fs)={\rm PD}^*_{\alpha,\theta}(d\fs)=\bE(\sigma_1^\theta;\sigma_1^{-1}\Delta\sigma_{[0,1]}\in d\fs)=\int_{(0,\infty)}x^\theta\bP(x^{-1}\Delta\sigma_{[0,1]}\in d\fs|\sigma_1=x)\bP(\sigma_1\in dx),\vspace{-0.1cm}$$
     where $(\sigma_t,t\ge 0)$ is a stable subordinator with Laplace transform $\bE(e^{-\lambda\sigma_t})=e^{-t\lambda^\alpha}$ and $\Delta\sigma_{[0,1]}$ the decreasing rearrangement of its jumps $\Delta\sigma_t=\sigma_t-\sigma_{t-}$, $t\in[0,1]$, see \cite{HPW,MPW} for details. 
\hfill $\square$ \end{ex}
There are simpler growth rules, which in general lead to models that are not fully exchangeable. Before we present these simpler growth rules, we mention a large class of models that retain a weak form of exchangeability and for which scaling limits have been obtained. 

\begin{ex}[Restricted exchangeable models \cite{CW}]\label{ex2b}\rm Let us define \em restricted exchangeable \em \linebreak[4] dislocation measures by their integral representation, referring to \cite{CW} for a full discussion:\vspace{-0.1cm}
    $$\kappa=c_1\delta_{\varepsilon^{(1)}}+\sum_{j\ge 1}\left(c_j\delta_{\varepsilon^{(j+1)}}+k_j\delta_{\omega^{[j]}}+\int_{\cS^\downarrow}\kappa_\fs(\,\cdot\,\cap\cP^{[j],\{j+1\}})\nu_j(d\fs)\right)\vspace{-0.2cm}$$
    where  
   $c_j
   \ge 0$, $k_j
   \ge 0$, $\omega^{[j]}=([j],\{j+1\},\{j+2\},\ldots)$, and 
    $\nu_j$ 
    is a measure on $\cS^\downarrow$ satisfying\vspace{-0.1cm}
	$$\nu_j(\{(1,0,\ldots),(0,0,\ldots)\})=0\quad\mbox{and}\quad\int_{\cS^\downarrow}\Bigg(s_01_{\{j=1\}}+\sum_{i\ge 1}s_i^j(1-s_i)\Bigg)\nu_j(d\fs)<\infty,\qquad j\ge 1.\vspace{-0.2cm}$$
  This includes all exchangeable dislocation measures, for $c_j\!=\!c$, $k_j\!=\!\nu(\{(0,\ldots)\})$, $\nu_j\!=\!\nu\!-\!k_j\delta_{(0,\ldots)}$. The splitting rules $(p_n,n\ge 2)$ associated with restricted exchangeable dislocation measures $\kappa$ give rise to the consistent restricted exchangeable 
  labelled Markov branching trees of \cite{CW}.

Consider the case where $c_j=k_j=0$ and $\lambda_n=\kappa(\cP\setminus\cP^{[n]})=n^\gamma\ell(n)$ for some slowly varying function $\ell$. The push-forward of $\kappa$ under $\Gamma\mapsto|\Gamma|^\downarrow$ is given by
$$\nu(d\fs)=\sum_{j\ge 1}\Bigg(\sum_{i\ge 1}s_i^j(1-s_i)\Bigg)\nu_j(d\fs).$$
Suppose that for each $j\geq 1$,  $\nu_j$ has its support in $\cS^\downarrow_1\setminus\{(1,0,\ldots)\}$ and that  $\int_{\cS^\downarrow}(1-s_1)\nu(d\fs)<\infty$.  Assuming further that $\nu_j=\nu_m$ for all $j\ge m$ for some $m\ge 1$, as in \cite[Theorem 7]{CW} where scaling limits were established
for convergence in probability, we deduce that \eqref{eq_tree1} holds for   
$$\kappa=\kappa_{\nu_m}-\kappa_{\nu_m}(\,\cdot\,\cap(\cP\setminus\cP^{[m]}))+\sum_{j=1}^{m-1}\int_{\cS^\downarrow}\kappa_\fs(\cdot\cap\cP^{[j],\{j+1\}})\nu_j(d\fs),$$
from the exchangeable case and by dominated convergence, because on the RHS only the measure $\kappa_{\nu_m}=\int_{\cS^\downarrow}\kappa_\fs(\cdot)\nu_m(d\fs)$ is infinite.
\hfill $\square$ \end{ex}   
One of the early families of regenerative tree growth processes to be studied was Ford's alpha-model. It has also been a main driver for much of the literature on scaling limits of Markov branching trees, both for general models and for further models with special structure. 

\begin{ex}[Ford's alpha-model \cite{For-05}]\rm
This family is parametrized by $\alpha\in [0,1]$ as follows.  For each edge $e$ of $\cT_n$, give $e$ weight $\alpha$ if both of its vertices are internal and weight $1-\alpha$ if one of its vertices is a leaf.  Choose an edge with probability proportional to its weight and attach $n+1$ to a new branch point between the two vertices of the selected edge.  From this description it is easy to check that $(\cT_n,n\geq 1)$ is a family of binary trees that forms a regenerative tree growth process.  Moreover, for $\pi=(B_1,B_2)$ we have
\[g_n(\pi,0)=\frac{\alpha}{n-\alpha}\qquad\mbox{and}\qquad g_n(\pi,i)=\frac{\#B_i-\alpha}{n-\alpha},\quad i\in\{1,2\}.\]
This model was introduced in \cite{For-05} as a model on cladograms that interpolates between the Yule model ($\alpha =0$), the uniform model ($\alpha =1/2$), and the comb ($\alpha =1$).  
\end{ex} 
The alpha-model is a restricted exchangeable model of binary trees that admits (at least) two natural extensions.  The alpha-gamma model, which is restricted exchangeable but not binary, and the alpha-theta model, which is binary but not, in general, restricted exchangeable.  The details of these models are our next two examples.

\begin{ex}[Alpha-gamma model \cite{CFW}]\label{ex1b}\rm
 For $0\le\gamma\le\alpha\le 1$ and $\pi=(B_1,\ldots,B_k)\in\cP_n$, let
	$$g_n(\pi,0)=\frac{\gamma}{n-\alpha},\qquad g_n(\pi,i)=\frac{\#B_i-\alpha}{n-\alpha},\ \ i\in[k],\qquad  g_n(\pi,k+1)=\frac{(k-1)\alpha-\gamma}{n-\alpha}.$$ The $\kappa$-measures are restricted exchangeable with\vspace{-0.1cm}
     $$c_j=k_j=0,\quad \mbox{and}\quad\nu_1(d\fs)=(1-\alpha){\rm PD}^*_{\alpha,-\alpha-\gamma}(d\fs),\quad\nu_j(d\fs)=\gamma{\rm PD}^*_{\alpha,-\alpha-\gamma}(d\fs),\ \ j\ge 2,\vspace{-0.1cm}$$
   if $0\!<\!\gamma\!<\!\alpha\!<\!1$, see \cite{CW}.  The convergence results of Example \ref{ex2b} include this as a special case.  
\hfill $\square$ \end{ex}

\begin{ex}[Alpha-theta model \cite{PW09}]\label{ex1a}\rm
  For $0\le\alpha\le 1$, $\theta\ge 0$ and $\pi=(B_1,B_2)\in\cP_n$, let
    $$g_n(\pi,0)=\frac{\alpha}{n-1+\theta},\quad g_n(\pi,1)=\frac{\#B_1-1+\theta}{n-1+\theta},\quad g_n(\pi,2)=\frac{\#B_2-\alpha}{n-1+\theta},\quad g_n(\pi,k)=0, \ k\ge 3.$$

This model is not restricted exchangeable except in the case $\theta=1-\alpha$, where the model reduces to Ford's alpha model.  Moreover, the dislocation measure for the alpha-theta model has not previously appeared in the literature.  To describe it, we introduce an \em ordered paintbox \em $\widetilde{\kappa}_{(u,1-u)}$, $0<u<1$, as the distribution of $\Pi=(\{i\ge 1\colon R_i=1\},\linebreak[1]\{i\ge 1\colon R_i=2\})$ where $R_1=1$ and the
$R_i$, $i\ge 2$, are independent random variables with $\bP(R_i\!=\!1)=u=1-\bP(R_i\!=\!2)$. 

For $0<\alpha<1$
and $\theta>0$, the $\kappa$-measure of the alpha-theta model is now given by\vspace{-0.1cm}
  $$\kappa=\alpha\widetilde{\kappa}^{\rm beta}_{\theta,-\alpha}(\,\cdot\,\cap\cP^{[2]})+\theta\widetilde{\kappa}^{\rm beta}_{\theta,-\alpha}(\,\cdot\,\cap\cP^{\{1\},\{2\}}),\ \mbox{where }\ \widetilde{\kappa}^{\rm beta}_{\theta,-\alpha}\!=\!\int_0^1\!\widetilde{\kappa}_{(u,1-u)}(\,\cdot\,)u^{\theta-1}(1-u)^{-\alpha-1}du.\vspace{-0.1cm}$$
To see this, note that from (\ref{splrules})\vspace{-0.1cm}
$$p_n(B_1,B_2)=
\left(\alpha1_{\{2\in B_1\}}+\theta1_{\{2\in B_2\}}\right)\frac{\Gamma(\#B_1-1+\theta)\Gamma(\#B_2-\alpha)}{\Gamma(n-1+\theta)\Gamma(1-\alpha)},\vspace{-0.1cm}$$
for $(B_1,B_2)\in\cP_n$ with $\#B_1\ge 1$ and $\#B_2\ge 1$, $n\ge 2$. The result now follows from (\ref{pkappa}) and the fact that $\widetilde{\kappa}_{\theta,-\alpha}^{\rm beta}(\cP^{B_1,B_2})\!=\!\int_0^1 u^{\theta-1+\#B_1-1}(1\!-\!u)^{-\alpha-1+\#B_2}du$ is a beta integral.    
    
Scaling limits for these trees were established in \cite{HM10} using criteria that are directly equivalent to what appears in our paper as Condition (i) in Theorem \ref{prop7} below.
  However, we can now give a shorter argument.  We have (similar to \cite{HM10} for exchangeable
  paintboxes)
  \[\int_\cP\left(1-|\Gamma^{[n]}|^\downarrow_1 \right)\widetilde{\kappa}_{(u,1-u)}(d\Gamma)\le \int_\cP\left(1-|\Gamma_1^{[n]}|\right)\widetilde{\kappa}_{(u,1-u)}(d\Gamma)=(1-u)\left(1-\frac{1}{n}\right) \]
But then \eqref{eq_tree1} follows for $\widetilde{\kappa}^{\rm beta}_{\theta,-\alpha}$ and for $\kappa$, which is bounded by a multiple of $\widetilde{\kappa}^{\rm beta}_{\theta,-\alpha}$, by dominated convergence. Thus Theorem \ref{thm_tree1} applies.
\hfill $\square$ \end{ex}
For the alpha-gamma and alpha-theta models, the regenerative property was shown
in \cite[Proposition 11]{PW09} and \cite[Proposition 8]{CW}, respectively. They both contain as special case for $\alpha=1/2$ and, respectively, $\theta=1/2$ and $\gamma=1/2$, the exchangeable uniform model on binary trees, related to Aldous's Brownian Continuum Random Tree \cite{Ald-91}.   Ford's binary alpha-model \cite{For-05} is also included in both examples.  

Aldous's binary beta model \cite{Ald-93} is included in the alpha-theta model for $\theta=-2\alpha$. Both the Poisson-Dirichlet model and the alpha-gamma model contain as special 
cases for $\alpha\in[1/2,1)$ and, respectively, $\theta=-1$ and $\gamma=1-\alpha$, the exchangeable model related to the stable Continuum Random Tree \cite{Du-LG,Mar-08,Mie-03}.

All the examples we have given so far satisfy  
the hypotheses of Theorem \ref{thm_tree1}. 
In fact, in these examples both tree convergence and residual mass process convergence were previously known to hold.  The
exchangeable case is \cite[Proposition 7]{HMPW}, the particle labelled 1 in the restricted exchangeable case is \cite[Proposition 28]{CW} and particle labelled 1 in the alpha-theta model is \cite[Proposition 6(iv)]{PW09}. Let us provide some very different examples that show what can go wrong.  We first give an example where tree convergence fails.

\begin{ex}\label{x1}\rm Some of the most elementary non-trivial dislocation measures are of the form
  $$\kappa=\sum_{j\ge 2}(\lambda_{j}-\lambda_{j-1})\delta_{\Gamma(j)}\qquad\mbox{for some $\Gamma(j)\in\cP^{[j-1],\{j\}}$, $j\ge 2$.}$$
  To ensure $\lambda_n\sim n^\gamma$ for some $\gamma\in(0,1)$, let $\lambda_n-\lambda_{n-1}=\gamma n^{\gamma-1}$.   
  For simplicity, we take 
  $\Gamma(j)$ binary with asymptotic frequencies $(x^{(j)},1-x^{(j)})$, where $x^{(j)}=1-1/j$. This implies
  $$\int_{\cS^\downarrow}(1-s_1)\nu(d\fs)=\int_{\cP}(1-|\Gamma|^\downarrow_1)\kappa(d\Gamma)=\sum_{j\ge 2}(1-x^{(j)})(\lambda_j-\lambda_{j-1})=\gamma\sum_{j\ge 2}j^{\gamma-2}<\infty,$$
  with $\nu$ as push-forward of $\kappa$. Consider a $\gamma$-self-similar tree $\cT_{\gamma,\nu}$ with dislocation measure $\nu$. 
  We explore two examples illustrating the validity/violation of 
  \eqref{eq_tree1}, which now reads  $$\int_\cP\left(|\Gamma^{[n]}|^\downarrow_1 -\left|\Gamma\right|^\downarrow_1\right)\kappa(d\Gamma)=\sum_{j\ge 2}\left(|\Gamma(j)^{[n]}|^\downarrow_1-x^{(j)}\right)(\lambda_j-\lambda_{j-1})\rightarrow 0,\quad\mbox{as $n\rightarrow\infty$.}$$  
 
 \begin{enumerate}\item[(a)] For $j\ge 2$ and $x^{(j)}=1-1/j\in(0,1)$, we construct $\Gamma(j)$ as a sequence $(\Gamma(j)^{[n]},n\ge j)$, starting from 
    $\Gamma(j)^{[j]}=([j-1],\{j\})$, and using Step $A_{x^{(j)}}$ inductively for $n\ge j$, where for $0\leq x\leq 1$:
  \begin{itemize}\item Step $A_x$: Given $\Gamma^{[n]}$, if $|\Gamma_1^{[n]}|>x$, set $\Gamma_1^{[n+1]}=\Gamma_1^{[n]}$, otherwise set 
      $\Gamma_1^{[n+1]}=\Gamma_1^{[n]}\cup\{n+1\}$. 
  \end{itemize}
  The purpose of Step $A_x$ is to change the relative frequency towards $x$. For $x=x^{(j)}$ and $\Gamma=\Gamma(j)$, we get $|\Gamma(j)_1|=|\Gamma(j)|^\downarrow_1=x^{(j)}$ 
  and 
$\left||\Gamma(j)_1^{[n]}|-x^{(j)}\right|\le 1-x^{(j)}$ for all $n\ge 1$, $j\ge 1$, equality for $j\ge n$ and strict inequality for $j<n$, since $1/n<1/j=1-x^{(j)}$. By the Dominated Convergence Theorem, \eqref{eq_tree1} is satisfied.\pagebreak[2]
  \item[(b)] For \eqref{eq_tree1} to fail, first let $\Gamma(j)_1^{[n]}$ approach frequency $1/2$, applying Step $A_{1/2}$ for $n<a_j$, so that $|\Gamma(j)_1^{[2j]}|=1/2$ and $|\Gamma(j)_1^{[n]}|\approx 1/2$ for $n\in[2j,a_j]$. Choose $(a_j)$ increasing with $2\ge\sum_{i\ge 2\colon n\in[2i,a_i]}(\lambda_i-\lambda_{i-1})>1$ and apply Step $A_{x^{(j)}}$ for $n\ge a_j$. Then 
    we will have $|\Gamma(j)_1|=x^{(j)}$ for all $j\ge 2$, but for all $n$ sufficiently large,
    $$\sum_{j\ge 2}\left(|\Gamma(j)_1^{[n]}|-x^{(j)}\right)(\lambda_j-\lambda_{j-1})\le-\frac{1}{3}\sum_{j\ge 2\colon n\in[2j,a_j]}(\lambda_j-\lambda_{j-1})<-\frac{1}{3}<0.$$
    Intuitively, the approximating trees have too many even branchpoints splitting into two equal-sized subtrees making trees wide and small in height, 
    while the proposed limiting distribution produces uneven branch points leading to thin and high trees with higher probability. Gromov-Hausdorff convergence fails, if total heights do not converge \cite{EPW}. \hfill $\square$ 
  \end{enumerate}
\end{ex}
In our next example, we show that the hypotheses of Theorem \ref{thm_mass1} are strictly stronger than the hypotheses of Theorem \ref{thm_tree1}.  

\begin{ex}\label{ex15}\rm In the general setting of Example \ref{x1}, consider $\Gamma(j)_1^{[n]}$ that first approaches (the wrong!) frequency
    $1-x^{(j)}$, applying Step $A_{1-x^{(j)}}$ for $n<a_j$, so that 
    $|\Gamma(j)_1^{[n]}|\approx 1-x^{(j)}$ for $n\in[j/(1-x^{(j)}),a_j]$. Then we apply Step $A_{x^{(j)}}$ for $n\ge a_j$ to achieve $|\Gamma(j)_1|=x^{(j)}$. We call these partitions ``evil''. If we did this for all $j\ge 2$, too many 
    partitions would have intermediate frequencies around 1/2 when restricted to $[n]$ and tree convergence may fail. Note that while at $1-x^{(j)}$, the block
    not containing 1 has frequency $x^{(j)}$ and is the larger block size that appears in the tree convergence criterion, while frequency $1-x^{(j)}$ is relevant for the residual mass process.

    To control the influence of partitions at intermediate frequencies, we also consider ``good'' partitions from Example \ref{x1}(a). The following strategy gives the right mix of ``good'' and
    ``evil'': 
    \begin{enumerate}\item[1.] For $j=2$ and $j=3$, start with two evil partitions $\Gamma(2)$ and 
        $\Gamma(3)$, with $|\Gamma(3)_1^{[\ell]}|\approx 1-x^{(3)}$ for $\ell=3/(1-x^{(3)})$, but
        leave $a_2$ and $a_3$ to be specified. Take good partitions $\Gamma(4),\ldots,\Gamma(\ell)$. 
        Also recall from the general setting that $\lambda_3-\lambda_2=\gamma 3^{\gamma-1}>0$. To
        proceed inductively, let $m=1$, $E_1=\{2,3\}$, $j_1=\ell+1$, and proceed to step 2.\vspace{-0.1cm}
      \item[2.] Given $(m,E_m,j_m)$, release the smallest evil partition $e_m=\min E_m$ by setting 
        $a_{e_m}=j_m$. Start evil partitions $\Gamma(j_m),\ldots,\Gamma(k_m)$ up to 
         $k_m\!=\!\inf\{j\!\ge\! j_m\colon\lambda_j\!-\!\lambda_{j_m}\ge\lambda_{e_m}\!\!-\!\lambda_{e_m\!-\!1}\}$. Let\vspace{-0.2cm}
        $$\ell_m=\inf\left\{j\ge k_m\colon|\Gamma(e_m)_1^{[j]}|\approx x^{(e_m)}\mbox{ and }|\Gamma(i)_1^{[j]}|\approx 1-x^{(i)},j_m\le i\le k_m\right\},\vspace{-0.2cm}
        $$
        and take good partitions $\Gamma(k_m+1),\ldots,\Gamma(\ell_m)$. Now set $E_{m+1}=(E_m\setminus\{e_m\})\cup\{j_m,\ldots,k_m\}$, $j_{m+1}=\ell_m+1$ and repeat step 2. for $(m+1,E_{m+1},j_{m+1})$.
    \end{enumerate}
    Now $|\Gamma(j)_1|=x^{(j)}$ for all $j\ge 2$ since $a_{e_m}<\infty$ for all evil partitions $e_m$. 
    The criterion \eqref{eq_tree1} of Theorem \ref{thm_tree1} for tree convergence holds, because the good
    partitions and the evil partitions that are
    either at frequency $x^{(j)}$ or $1-x^{(j)}$ give convergence as in Example \ref{x1}(a), while the evil 
    partitions at intermediate frequencies have total weight
    $w_m=(\lambda_{e_m}-\lambda_{e_m-1})+(\lambda_{k_m}-\lambda_{j_m-1})\rightarrow 0$ as $m\rightarrow\infty$,
    so their contribution vanishes as $m\rightarrow\infty$.

    The criterion \eqref{eq_mass1} of Theorem \ref{thm_mass1} for residual mass process convergence is not satisfied,
    because for every $n\ge 3/(1-x^{(3)})$, there are evil partitions of weight at least 
    $\lambda_3-\lambda_2$ which have a frequency $|\Gamma(j)^{[n]}_1|\approx 1-x^{(j)}$ that is smaller by more than 1/4 than their limit frequency $x^{(j)}$,
    since $x^{(j)}-(1-x^{(j)})>1/4$ for all $j\ge 3$, and this cannot be offset by partitions that exceed their limit frequencies, by the argument in Example \ref{x1}(a).
\hfill $\square$ \end{ex}

\section{Background}\label{sec_background}

In this section we present the background information needed to understand the statements of our results on scaling limits of random trees.  Since the proofs of our results do not require any technical details about the constructions in this section we keep the discussion light and heuristic at times, referring to the existing literature for details. 

\subsection{Trees as metric measure spaces}
The trees under discussion in this paper can naturally be considered as metric spaces with the graph metric.  That is, the distance between two vertices is the number of edges on the path connecting them.  Let $(T,d, \textrm{root})$ be a tree equipped with the graph metric.  For $a>0$, we define $at$ to be the metric space $(T,ad,\textrm{root})$, i.e. the metric is scaled by $a$.  Moreover, the trees we are dealing with are rooted so we consider $(T,d,\textrm{root})$ as a pointed metric space with the root as the point.  Additionally, we let $\mu_T$ be the uniform probability measure on the leaves of $T$.  If we have a random tree $\cT$, this gives rise to a random pointed metric measure space $(\cT,d,\textrm{root},\mu_{\cT})$.  For this last statement to be made rigorous, it is clear that we need to put a topology on pointed metric measure spaces.  This is hard to do in general, but note that the pointed metric measure spaces that come from the trees we are discussing are compact and this simplifies matters.

Let $\M$ be the set of equivalence classes of compact pointed metric measure spaces (equivalence here being up to point and measure preserving isometry).  We endow $\M$ with the pointed Gromov-Hausdorff-Prokhorov metric (see \cite{HM10}).  Fix $(X,d,\rho,\mu), (X',d',\rho,\mu') \in \M$ and define\vspace{-0.1cm}
\[ d_{\textrm{GHP}}(X,X') = \inf_{(M,\delta)} \inf_{\phi:X\to M \atop \phi' :X'\to M}\left[ \delta(\phi(\rho),\phi'(\rho')) \vee \delta_H(\phi(X),\phi'(X'))\vee \delta_P(\phi_*\mu,\phi'_*\mu')\right],\vspace{-0.2cm}\]
where the first infimum is over metric spaces $(M,\delta)$, the second infimum if over isometric embeddings $\phi$ and $\phi'$ of $X$ and $X'$ into $M$, $\delta_H$ is the Hausdorff distance on compact subsets of $M$, and $\delta_P(\phi_*\mu,\phi'_*\mu')$ is the Prokhorov distance between the push-forward $\phi_*\mu$ of $\mu$ by $\phi$ and the push-forward $\phi'_*\mu'$ of $\mu'$ by $\phi'$. It is worth noting briefly that the definitions of $\M$ and $d_{\textrm{GHP}}$ as just given do not make formal sense in Zermelo-Fraenkel set theory with the axiom of choice (ZFC); one might just as well try metrizing the set of all sets. \pagebreak[4] Nonetheless, it is not hard to formalize the heuristic definitions we have given.  For example, one can use the fact that every separable metric space can be isometrically embedded in $\ell_\infty$ to find an honest set $\overline{\mathcal{M}}_w$ in ZFC such that every compact pointed metric measure space is isometric, by point and measure preserving isometry, to exactly one element of $\overline{\mathcal{M}}_w$ and then do everything internally in this set.

\begin{prop}[Proposition 1 in \cite{HM10}, see also \cite{EW06, GPW09, G99, M09}] 
The space $(\M, d_{\textrm{GHP}})$ is a complete separable metric space.
\end{prop}
Scaling limits of discrete trees are elements of $\M$ that are tree-like metric spaces.  An $\bR$-tree is a complete metric space $(T,d)$ with the following properties:\vspace{-0.1cm}
\begin{itemize}
\item For $v,w\in T$, there exists a unique isometry $\phi_{v,w}\colon[0,d(v,w)]\to T$ with $\phi_{v,w}(0)=v$ and $\phi_{v,w}(d(v,w))=w$.\vspace{-0.2cm}
\item For every continuous injective function $c\colon[0,1]\to T$ such that $c(0)=v$ and $c(1)=w$, we have $c([0,1]) = \phi_{v,w}([0,d(v,w)])$.\vspace{-0.1cm}
\end{itemize}
If $(T,d)$ is a compact $\bR$-tree, every choice of root $\rho\in T$ and probability measure $\mu$ on $T$ yields an element $(T,d,\rho,\mu)$ of $\M$.  With this choice of root also comes a height function $\Ht(v) = d(v,\rho)$.  The leaves of $T$ can then be defined as the points $v\in T$ such that $v$ is not in $[[\rho,w[[ \deq \phi_{\rho,w}([0,\Ht(w)))$ for any $w\in T$.  The set of leaves is denoted $\mathcal{L}(T)$.  

\begin{defn}\rm
A \textit{continuum tree} is an $\bR$-tree $(T,d,\rho,\mu)$ with a choice of root and probability measure such that $\mu$ is non-atomic, $\mu(\mathcal{L}(T))=1$, and for every non-leaf vertex $w$, $\mu(\{v\in T\!\!:[[\rho,v]]\cap [[\rho,w]] = [[\rho,w]]\})>0$.
\end{defn}
A continuum random tree (CRT) is an $(\M,d_{\textrm{GHP}})$-valued random variable that is almost surely a continuum tree.  The continuum random trees we will be interested in are those associated with self-similar mass fragmentation processes.    
\vspace{-0.2cm}

\subsection{Self-similar mass fragmentations}\label{subsec_ft}
We are now prepared to introduce self-similar mass fragmentations and their genealogical trees.  Suppose $\gamma>0$ and let $\nu$ be a $\sigma$-finite measure on $\cS^\downarrow$ such that $\nu(\{(1,0,0,\dots)\})=0$ and $\int_{\cS^\downarrow} (1-s_1)\nu(d\mathbf{s}) < \infty$ and $\nu(\sum_i s_i <1)=0$.  Heuristically, a self-similar mass fragmentation with characteristics $(\gamma,\nu)$ is an $\cS^\downarrow$-valued Markov process $(F(t),t\geq 0)$ such that $F(0)=(1,0,0,\dots)$ and such that a block of size $x$ splits into blocks $x\mathbf{s} = (xs_1,xs_2,\dots)$ at rate $x^{-\gamma} \nu(d\mathbf{s})$.  A rigorous construction of such processes can be found in \cite{Ber-book}, though we remark that there is a slight difference in notation: our index $\gamma$ of self-similarity corresponds to the index $-\gamma$ in \cite{Ber-book}.  The idea of the genealogical tree of a self-similar mass fragmentation is to construct an $\bR$-tree that keeps track of the sizes of the blocks of the fragmentation process as time progresses.

For a continuum tree $(T,\mu)$ and $t\geq 0$, let $T_1(t),T_2(t),\dots$ be the tree components of $\{v\in T\!: \mathrm{ht}(v)>t\}$, ranked in decreasing order of $\mu$-mass (breaking ties uniformly).  We call a continuum random tree $(\cT,\mu)$ $\gamma$-self-similar if for every $t\geq 0$, conditionally on $(\mu(\cT_i(t)),i\geq 1)$, $(\cT_i(t),i\geq 1)$ has the same law as $(\mu(\cT_i(t))^{\gamma}\cT^{(i)},i\geq 1)$ where the $\cT^{(i)}$, $i\ge 1$, are independent copies of $\cT$. 

The following summarizes the parts of Theorem 1 and Lemma 5 in \cite{HM04} that we will need.

\begin{thm}\label{theorem existence of trees}
Let $F$ be a $(\gamma,\nu)$-self-similar fragmentation with $\gamma>0$ and $\nu$ as above.  There exists a $\gamma$-self-similar CRT $(\cT_{\gamma ,\nu},\mu_{\gamma,\nu})$ such that, writing $F'(t)$ for the decreasing sequence of masses of the connected components of $\{v\in \cT_{\gamma ,\nu}\colon{\rm ht}(v)>t\}$, the process $(F'(t),t\geq 0)$ has the same law as $F$.  Furthermore, $\cT_{\gamma,\nu}$ is a.s.\ compact.
\end{thm}  
The Brownian continuum random tree introduced by Aldous \cite{Aldo93} as the scaling limit of conditioned Galton-Watson trees is an example of a self-similar fragmentation tree.  

\begin{defn}\rm The Brownian CRT is the $1/2$-self-similar random tree with dislocation measure $\nu$ given by\vspace{-0.4cm}
\[\int_{\cS^\downarrow} f(\mathbf{s}) \nu(d\mathbf{s}) = \int_{1/2}^1\sqrt{\frac{2}{\pi s_1^3(1-s_1)^3}}  f(s_1,1-s_1,0,0,\dots)ds_1.\]
\end{defn}\pagebreak
One of our main tools will be the general theory of scaling limits of unordered Markov branching trees.  In particular, we make use of the following theorem.

\begin{thm}[Theorem 5 in \cite{HM10}] \label{thm_HM10}
Let $(\cT^\circ_n,n\geq 1)$ be a Markov branching model based on $(p_n^\circ,n\geq 2)$ as in Section \ref{subsec M-branch}.  Suppose that there is a characteristic pair $(\gamma,\nu)$ with $\gamma>0$, and $\nu$ satisfying the conditions at the start of Section \ref{subsec_ft} as well as a function $\ell\colon(0,\infty)\to (0,\infty)$, slowly varying at $\infty$ such that, in the sense of weak convergence of finite measures on $\cS^\downarrow$, we have
\begin{equation}\label{eq_HMcond} n^\gamma \ell(n)(1-s_1)\bar p_n^\circ(d\mathbf{s}) \rightarrow (1-s_1)\nu(d\mathbf{s}),\end{equation}
where $\bar p_n^\circ$ is the push-forward of the measure on $\cP_n^\circ$ with probability function $p_n^\circ$ onto $\cS^\downarrow$ by the map
\[ (n_1,\dots, n_p) \mapsto \left(\sum_{i=1}^p n_i\right)^{-1} (n_1,\dots, n_p,0,0,\dots).\] 
If we view $\cT^\circ_n$ as a random element of $\M$ with the graph distance and the uniform probability measure its leaves,
then we have the convergence in distribution
\[ \frac{1}{n^\gamma \ell(n)} \cT^\circ_n \rightarrow \cT_{\gamma,\nu},\]
with respect to the rooted Gromov-Hausdorff-Prokhorov topology.
\end{thm}

\section{Scaling limits of regenerative tree growth processes}\label{subsecHMtrees}

While every dislocation measure $\kappa$ on $\cP$ gives rise to a regenerative tree growth process, not every such process has a scaling limit. Examples without scaling limit include the $(0,\theta)$-tree growth process studied in \cite[Proposition 13]{PW09}, where the growth is logarithmic and the branching structure degenerates under logarithmic scaling. In view of Proposition \ref{prop_mb}, it makes sense to try interpreting the hypotheses of Theorem \ref{thm_HM10} in terms of $\kappa$.  In particular, let us examine the LHS of \eqref{eq_HMcond}.  From Proposition \ref{prop_mb} we see that for a bounded continuous function $f$ on $\cS^\downarrow$, 
\[\begin{split} &\hspace{-0.3cm}\int_{\cS^\downarrow}f(\fs)n^\gamma\ell(n)(1-s_1)\bar p_n^\circ(d\fs) \\
	& = \sum_{n_1\ge\cdots\ge n_k:n_1+\cdots+n_k=n}n^\gamma\ell(n)p_n^\circ(n_1,\ldots,n_k)\left(1-\frac{n_1}{n}\right)f\left(\frac{n_1}{n},\ldots,\frac{n_k}{n},0,\ldots\right) \\
&=\frac{n^\gamma\ell(n)}{\lambda_n} \sum_{\pi\in\cP_n\setminus\{1_{[n]}\}}\kappa(\cP^\pi)\left(1-\frac{(\#\pi)^\downarrow_1}{n}\right)f\left(\frac{(\#\pi)^\downarrow}{n}\right)\\ &=\frac{n^\gamma\ell(n)}{\lambda_n} \int_{\cP}\left(1-|\Gamma^{[n]}|^\downarrow_1 \right)f\left(|\Gamma^{[n]}|^\downarrow \right)\kappa(d\Gamma)
,\end{split}\]
where now we write $(\#\pi)^\downarrow$ for the decreasing rearrangement of the block sizes of $\pi$, with an infinite string of zeros appended (whereas in our previous usage $(\#\pi)^\downarrow$ was a finite vector).  Given this expression and the convergence \eqref{eq_HMcond} we need to establish, natural assumptions on $\kappa$ become that $\lambda_n = n^\gamma \ell(n)$ for some $\gamma>0$ and $\ell(n)$ slowly varying at $\infty$ and that
\[ \lim_{n\to\infty} \int_{\cP}\left(1-|\Gamma^{[n]}|^\downarrow_1 \right)f\left(|\Gamma^{[n]}|^\downarrow \right)\kappa(d\Gamma) = \int_{\cP}\left(1-|\Gamma|^\downarrow_1 \right)f\left(|\Gamma|^\downarrow \right)\kappa(d\Gamma). \]
Of course, for this last equation to have hope of holding, we must assume that $\kappa$-a.e.\ $\Gamma\in\cP$ has asymptotic ranked frequencies.  This holds for exchangeable and restricted exchangeable $\kappa$, and when $\kappa$ is partially exchangeable in the sense of \cite{Pit-95}.  Let us, though, clarify the relationship between the existence of asymptotic frequencies and the existence of asymptotic ranked frequencies.

\begin{lm}\label{lm7} Existence of $|\Gamma|^\downarrow_i$ for \em all \em $i\ge 1$ holds if and only if 
  $(|\Gamma_i|,i\ge 1)$ exists as a uniform limit.\linebreak In this case, $(|\Gamma|^\downarrow_i,i\ge 1)$ is the decreasing rearrangement of $(|\Gamma_i|,i\ge 1)$, which we write as $|\Gamma|^\downarrow$.
\end{lm}
Lemma \ref{lm7} is inessential to the remainder of our results, but for completeness we include a proof in Appendix \ref{ap_lemma7}. Note that asymptotic (ranked) frequencies need not be in $\cS^\downarrow_1$ and that $|\Gamma_i|$ may vanish. 
The partition $\Gamma=(\{2^{i-1},\ldots,2^i-1\},i\ge 1)$ is an example where $|\Gamma_i|$, $i\ge 1$, exists, but $|\Gamma|^{\downarrow}_1$ does not. 

We can now give our main result on the existence of scaling limits of regenerative growth processes, which contains the statement of Theorem \ref{thm_tree1}.
\begin{thm}\label{prop7} Let $(\cT_n,n\geq 1)$ be a regenerative tree growth process associated with a dislocation measure $\kappa$.  Assume that $\kappa$-a.e. 
      $\Gamma\in\cP$ has proper asymptotic ranked frequencies in $\cS^\downarrow_1\setminus\{(1,0,\ldots)\}$, define
      $\nu$ to be the push-forward of $\kappa$ under $\Gamma\mapsto|\Gamma|^\downarrow$ and suppose  $\int_{\cS^\downarrow}(1-s_1)\nu(d\fs)<\infty$ and $\lambda_n=\kappa(\cP\setminus\cP^{[n]})=n^\gamma\ell(n)$. Then 
      the following are equivalent:
  \begin{enumerate}\item[\rm(i)] For all bounded continuous $f\colon\cS^\downarrow\rightarrow[0,\infty)$,\vspace{-0.1cm} 
	  $$\int_{\cP}\left(1-|\Gamma^{[n]}|^\downarrow_1 \right)f\left(|\Gamma^{[n]}|^\downarrow\right)\kappa(d\Gamma)
			\rightarrow\int_{\cS^\downarrow}(1-s_1)f(\fs)\nu(d\fs),\qquad\mbox{as $n\rightarrow\infty$}\textit{;}\vspace{-0.3cm}$$
    \item[\rm(ii)] $\displaystyle\int_{\cP}\left(|\Gamma^{[n]}|^\downarrow_1 -|\Gamma|^\downarrow_1\right)\kappa(d\Gamma)\rightarrow 0$, as $n\rightarrow\infty$;\vspace{-0.1cm}
    \item[\rm(iii)] $\displaystyle\int_{\cP}\left||\Gamma^{[n]}|^\downarrow_1 -|\Gamma|^\downarrow_1\right|\kappa(d\Gamma)\rightarrow 0$, as $n\rightarrow\infty$, i.e.\ the convergence $|\Gamma^{[n]}|^\downarrow_1 \rightarrow|\Gamma|^\downarrow_1$ holds in $L_1(\kappa)$.\vspace{-0.1cm}
  \end{enumerate}
  If condition {\rm(ii)} holds, then $\displaystyle\frac{\cT_n^\circ}{n^\gamma\ell(n)}\rightarrow\cT_{\gamma,\nu}$ in distribution, as $n\rightarrow\infty$, in the rooted GHP sense.
\end{thm}
\begin{pf} If (i) holds, we obtain (ii) as a rearrangement of the special case $f=1$. Now assume (ii). Let us prove (iii).   The main difficulty arises from the possibility that $\kappa(\cP^{[m]})=\infty$.  For all $m\ge1$ we have
  $$
                 \int_\cP\left||\Gamma^{[n]}|^\downarrow_1-|\Gamma|^\downarrow_1\right|\kappa(d\Gamma)=\!\int_{\cP^{[m]}}\left||\Gamma^{[n]}|^\downarrow_1-|\Gamma|^\downarrow_1\right|\kappa(d\Gamma)+\int_{\cP\setminus\cP^{[m]}}\left||\Gamma^{[n]}|^\downarrow_1-|\Gamma|^\downarrow_1\right|\kappa(d\Gamma).
  $$
Since $\kappa(\cP\setminus\!\cP^{[m]})\!<\!\infty$ and $\kappa$-a.e.\ $\Gamma$ has asymptotic ranked frequencies an application of dominated convergence shows that, for each fixed $m$, the second term vanishes as $n\rightarrow\infty$.  From the triangle inequality, we see that
\[\begin{split}
 \int_{\cP^{[m]}}\left||\Gamma^{[n]}|^\downarrow_1-|\Gamma|^\downarrow_1\right|\kappa(d\Gamma) & = \int_{\cP^{[m]}}\left|(1-|\Gamma^{[n]}|^\downarrow_1)-(1-|\Gamma|^\downarrow_1)\right|\kappa(d\Gamma)\\
 &\hspace{-0.5cm} \leq \int_{\cP^{[m]}} \left(1\!-\!|\Gamma^{[n]}|^\downarrow_1\right)\!\kappa(d\Gamma) + \int_{\cP^{[m]}}\!\left(1\!-\!|\Gamma|^\downarrow_1\right)\!\kappa(d\Gamma) \\
 &\hspace{-0.5cm} = \int_{\cP}\!\left(1\!-\!|\Gamma^{[n]}|^\downarrow_1\right)\!\kappa(d\Gamma)-\!\int_{\cP\setminus\cP^{[m]}}\!\!\left(1\!-\!|\Gamma^{[n]}|^\downarrow_1\right)\!\kappa(d\Gamma)+\!\int_{\cP^{[m]}}\!\left(1\!-\!|\Gamma|^\downarrow_1\right)\!\kappa(d\Gamma).
\end{split}\]
It follows from (ii) that 
\[\lim_{n\to\infty}  \int_{\cP}\!\left(1\!-\!|\Gamma^{[n]}|^\downarrow_1\right)\!\kappa(d\Gamma) =  \int_{\cP}\!\left(1\!-\!|\Gamma|^\downarrow_1\right)\!\kappa(d\Gamma)\]
and, since $\kappa(\cP\setminus\!\cP^{[m]})\!<\!\infty$, dominated convergence shows that
\[\lim_{n\to\infty} \int_{\cP\setminus\cP^{[m]}}\!\!\left(1\!-\!|\Gamma^{[n]}|^\downarrow_1\right)\!\kappa(d\Gamma) =\int_{\cP\setminus\cP^{[m]}}\!\!\left(1\!-\!|\Gamma|^\downarrow_1\right)\!\kappa(d\Gamma).\]
Consequently, for every $m\geq 1$ we have
\[\limsup_{n\to\infty} \int_{\cP}\left||\Gamma^{[n]}|^\downarrow_1-|\Gamma|^\downarrow_1\right|\kappa(d\Gamma) \leq \limsup_{n\to\infty}  \int_{\cP^{[m]}}\left||\Gamma^{[n]}|^\downarrow_1-|\Gamma|^\downarrow_1\right|\kappa(d\Gamma) \leq 2\int_{\cP^{[m]}}\!(1\!-\!|\Gamma|^\downarrow_1)\kappa(d\Gamma).\]
Since $\int_{\cP}(1-|\Gamma|^\downarrow_1)\kappa(d\Gamma)<\infty$, while $\bigcap_{m\ge 1}\cP^{[m]}=\{1_\bN\}$ and $\kappa(\{1_\bN\})=0$, the infimum of these bounds over $m\ge 1$ vanishes and (iii) follows.
  
Now assume (iii). If $f\colon\cS^\downarrow\rightarrow[0,\infty)$ is continuous and
   bounded, then
  \begin{eqnarray*}&&\hspace{-0.5cm}\left|\int_{\cP}\left(1-|\Gamma^{[n]}|^\downarrow_1 \right)f\left(|\Gamma^{[n]}|^\downarrow\right)\kappa(d\Gamma)
		       -\int_{\cP}(1-|\Gamma|^\downarrow_1)f(|\Gamma|^\downarrow)\kappa(d\Gamma)\right|\\
    &&\le\int_{\cP}\left||\Gamma^{[n]}|^\downarrow_1-|\Gamma|^\downarrow_1\right|\left|f(|\Gamma^{[n]}|^\downarrow)\right|\kappa(d\Gamma)+\int_{\cP}(1-|\Gamma|^\downarrow_1)\left|f(|\Gamma^{[n]}|^\downarrow)-f(|\Gamma|^\downarrow)\right|\kappa(d\Gamma),
  \end{eqnarray*}
  and (i) follows from (iii) by dominated convergence, since $\kappa$-a.e.\ $\Gamma$ has asymptotic ranked frequencies, since $f$ is bounded and continuous, and since $\nu$ is the push-forward of $\kappa$ under $\Gamma\mapsto|\Gamma|^\downarrow$.

  The last part follows from Haas and Miermont \cite[Theorem 5]{HM10}, which we formulated in Theorem \ref{thm_HM10} above.
\end{pf}

\section{Residual mass processes in regenerative tree growth processes}\label{secHMmcs}

Let $(\cT_n,n\ge 1)$ be a regenerative tree growth process and $(X_m^{(n)},m\ge 0)$ the associated residual mass processes of label 1 in $\cT_n$, $n\ge 1$, with transition probabilities 
$$\bP(X^{(n)}_1=k)=\sum_{\pi=(B_1,\ldots,B_k)\in\cP_n:\#B_1=k}p_n(\pi)=\frac{1}{\lambda_n}\kappa\left(\left\{\Gamma\in\cP\colon\#\Gamma_1^{[n]}=k\right\}\right),\quad 1\le k<n,$$
as identified in Proposition \ref{resmp}, with $\lambda_n=\kappa(\{\Gamma\in\cP\colon\Gamma_1^{[n]}\neq[n]\})$. The existence of a scaling limit for
trees $\cT_n^\circ$ as studied in Section \ref{subsecHMtrees} does not imply the existence of a scaling limit for associated residual mass processes $X^{(n)}$, in general (see Example \ref{ex15}). In this section, we study scaling limits $X_{\lfloor\lambda_n t\rfloor}^{(n)}/n\rightarrow X_t$, as $n\rightarrow\infty$.  Since, for fixed $n\ge 1$, $(X^{(n)}_m, m\geq 0)$ is a non-increasing Markov chain with $X^{(n)}_0=n$, we can make use of the general theory of self-similar scaling limits for such chains that was recently developed in \cite{HM09}. 

\begin{thm}[Theorems 1 and 2 in \cite{HM09}]\label{thm_HM09}
Let $p=(p_{ij}, 0\leq j\leq i)$ be a transition matrix, and for each $n\ge 1$ let $(Y^{(n)}_m, m\geq 0)$ be a Markov chain with transition matrix $p$ such that $Y^{(n)}_0=n$.  Define\vspace{-0.2cm}
\[p^*_n(dx) = \sum_{k=0}^n p_{n,k} \delta_{k/n}(dx).\]
Suppose that there exists a sequence $(a_n,n\geq 0)$ of the form $a_n=n^\gamma\ell(n)$ for some $\gamma>0$ and a slowly varying function $\ell$ as well as a non-zero finite measure $\mu$ on $[0,1]$ such that
\begin{equation}\label{eq_massh} a_n(1-x)p^*_n(dx) \rightarrow \mu(dx)\end{equation}
in the sense of weak convergence of finite measures on $[0,1]$.  Then we have the following convergence in distribution\vspace{-0.2cm}
\[ \left(\frac{Y^{(n)}_{\lfloor a_n t\rfloor}}{n}, \ t\geq 0\right) \rightarrow \left(X_t, \ t\geq 0\right)\]
in the Skorokhod sense, where $X$ is a self-similar Markov process and in Lamperti's representation {\rm(\ref{pssmp})}, we have
\[ \bE(e^{-s\xi_r})=\exp(-r\psi(s))\qquad\mbox{with}\quad\psi(s) = \int_{[0,1]} \frac{1-x^s}{1-x} \mu(dx).\]
Moreover, letting $A_n$ be the absorption time of $Y^{(n)}$ at $0$, the above convergence in distribution happens jointly with the convergence of $a_n^{-1}A_n$ to the absorption time at $0$ of the limiting process. 
\end{thm}\pagebreak[2]   
Proposition \ref{resmp} shows that the residual mass process of the leaf $\{1\}$ falls into the scope of this theorem with $p_{ij}=\bP(X^{(i)}_1=j)$.
Given a dislocation measure $\kappa$ with $\lambda_n=\kappa(\cP\setminus\cP^{[n]})=n^\gamma\ell(n)$ regularly varying, as $n\rightarrow\infty$, for some $\gamma > 0$, and taking $a_n=\lambda_n$, the LHS of condition (\ref{eq_massh}) for the residual mass process of leaf $\{1\}$ becomes\vspace{-0.1cm}
$$\sum_{\pi=(B_1,\ldots,B_k)\in\cP_n}n^\gamma\ell(n)p_n(\pi)\left(1-\frac{\#B_1}{n}\right)f\left(\frac{\#B_1}{n}\right)=\int_\cP\! \left(1-|\Gamma_1^{[n]}| \right)f\left(|\Gamma_1^{[n]}| \right)\kappa(d\Gamma).
$$
\begin{thm}\label{propmc} Let $(\cT_n,n\geq 1)$ be a regenerative tree growth process with dislocation measure $\kappa$ and $X^{(n)}$ the residual mass process of $\{1\}$ in $\cT_n$.  Assume that the first block $\Gamma_1$ of $\kappa$-a.e.\ $\Gamma\in\cP$ has an asymptotic frequency $|\Gamma_1|\in(0,1)$, define $\Lambda$ as push-forward of $\kappa$ under $\Gamma\mapsto-\log(|\Gamma_1|)$. If $\int_{(0,\infty)}(1-e^{-x})\Lambda(dx)<\infty$ and $\lambda_n=\kappa(\cP\setminus\cP^{[n]})=n^\gamma\ell(n)$, then the following are equivalent:
\begin{enumerate}\item[\rm(i)] For all bounded continuous $f\colon[0,1]\rightarrow [0,\infty)$,\vspace{-0.1cm}
	$$\int_\cP \left(1-|\Gamma_1^{[n]}| \right)f\left(|\Gamma_1^{[n]}| \right)\kappa(d\Gamma)\rightarrow\int_{(0,\infty)}f(e^{-y})(1-e^{-y})\Lambda(dy),\ \mbox{as $n\rightarrow\infty$;}\vspace{-0.3cm}
$$
  \item[\rm(ii)] $\displaystyle\int_{\cP}\left(|\Gamma_1^{[n]}| -|\Gamma_1|\right)\kappa(d\Gamma)\rightarrow 0$, as $n\rightarrow\infty$;\vspace{-0.1cm}
  \item[\rm(iii)] $\displaystyle\int_{\cP}\left||\Gamma_1^{[n]}| -|\Gamma_1|\right|\kappa(d\Gamma)\rightarrow 0$, as $n\rightarrow\infty$, i.e.\ the convergence $|\Gamma_1^{[n]}|\rightarrow|\Gamma_1|$ holds in $L_1(\kappa)$. 
  \end{enumerate}
  If condition {\rm(ii)} holds, $X^{(n)}_{\lfloor\lambda_n t\rfloor}/n\!\rightarrow\! X_t$ in distribution, as $n\!\rightarrow\!\infty$, in the Skorohod sense as functions of $t\ge 0$, where $X$ is a self-similar Markov process and  $\bE(e^{-s\xi_r})=\exp(-r\int_{(0,\infty)}(1\!-\!e^{-sy})\Lambda(dy))$ in Lamperti's representation {\rm(\ref{pssmp})}.  Moreover, letting $A_n$ be the absorption time of $X^{(n)}$ at $0$, the above convergence in distribution happens jointly with the convergence of $\lambda_n^{-1}A_n$ to the absorption time at $0$ of the limiting process.   

  If, in addition, $\kappa$-a.e.\ $\Gamma\in\cP$ has asymptotic ranked  
  frequencies, then \eqref{eq_tree1} holds and we have $\displaystyle\frac{\cT_n^\circ}{n^\gamma\ell(n)}\rightarrow\cT_{\gamma,\nu}$ in distribution, as $n\rightarrow\infty$, in the rooted GHP sense.
\end{thm}
\medskip

\noindent We note that the statement of Theorem \ref{thm_mass1} is contained in the statement of this theorem.

\begin{pf} The proof of the equivalences is the same as for Theorem \ref{prop7}, with $|\Gamma^{[n]}|^\downarrow_1$ and $|\Gamma|^\downarrow_1$ replaced by $|\Gamma_1^{[n]}|$ and $|\Gamma_1|$. Convergence of $X^{(n)}_{\lfloor\lambda_n t\rfloor}/n$  is now an application of \cite[Theorem 1]{HM09}, which we formulated as Theorem \ref{thm_HM09} above. Finally, if $\kappa$-a.e.\ $\Gamma\in\cP$ has asymptotic ranked frequencies, we first note that in the notation of Theorem \ref{prop7}\vspace{-0.1cm}
\begin{equation}\label{lambdanu}\int_{S^{\downarrow}}(1-s_1)\nu(d\fs)=\int_\cP(1-|\Gamma|^\downarrow_1)\kappa(d\Gamma)\le\int_\cP(1-|\Gamma_1|)\kappa(d\Gamma)=\int_{(0,\infty)}(1-e^{-x})\Lambda(dx)<\infty.\vspace{-0.1cm}\end{equation}
To apply Theorem \ref{prop7}, we verify condition (ii) of Theorem \ref{prop7}:\vspace{-0.1cm} 
  \begin{eqnarray} \int_\cP\left(|\Gamma^{[n]}|^\downarrow_1-|\Gamma|^\downarrow_1\right)\kappa(d\Gamma)\!&\!\!\!=\!\!\!&\!\!\int_{\{|\Gamma_1|=|\Gamma|^\downarrow_1\}}\left(|\Gamma_1^{[n]}| -|\Gamma_1|\right)\kappa(d\Gamma)
			+\int_{\{|\Gamma_1|\neq|\Gamma|^\downarrow_1\}}\left(|\Gamma^{[n]}|^\downarrow_1 -|\Gamma|^\downarrow_1\right)\kappa(d\Gamma)\nonumber\\
		&&\hspace{0.5cm}+\int_{\{|\Gamma_1|=|\Gamma|^\downarrow_1\}}\left(|\Gamma^{[n]}|^\downarrow_1-|\Gamma_1^{[n]}|\right)\kappa(d\Gamma)\label{split}\vspace{-0.1cm}
\end{eqnarray}
is a sum of three terms. The first term vanishes as $n\rightarrow\infty$ by (iii). The second term 
  vanishes as $n\rightarrow\infty$ since $\kappa(|\Gamma_1|\neq|\Gamma|^\downarrow_1)<\infty$: if $|\Gamma_1|\neq|\Gamma|^\downarrow_1$, then one of them must be less than 1/2,
  so $\kappa(|\Gamma_1|\neq|\Gamma|^\downarrow_1)\le\kappa(|\Gamma_1|\le 1/2)+\kappa(|\Gamma|^\downarrow_1\le 1/2)<\infty$, by (\ref{lambdanu}). The third term is non-negative, so that $\liminf_{n\rightarrow\infty}{\rm LHS}\ge 0$ in (\ref{split}). In\vspace{-0.1cm}
$$
                 \int_\cP\left(|\Gamma^{[n]}|^\downarrow_1-|\Gamma|^\downarrow_1\right)\kappa(d\Gamma)\le\!\int_{\cP^{[m]}}\left(1-|\Gamma_1|\right)\kappa(d\Gamma)+\int_{\cP\setminus\cP^{[m]}}\left(|\Gamma^{[n]}|^\downarrow_1-|\Gamma|^\downarrow_1\right)\kappa(d\Gamma)\vspace{-0.1cm}$$
  we can make the first term small by choosing $m$ large and 
  the second term vanishes as $n\rightarrow\infty$, for each fixed $m\ge 1$. Hence, $\limsup_{n\rightarrow\infty} {\rm LHS}\le 0$ and so $\lim_{n\rightarrow\infty} {\rm LHS}=0$. 
\end{pf}\pagebreak

\noindent Residual mass convergence and tree convergence are not equivalent. The last part of Theorem \ref{propmc} finds that under the conditions for residual mass convergence in this theorem, we just need to assume the existence of asymptotic ranked frequencies to also obtain tree convergence. Example \ref{ex15} demonstrates that residual mass convergence does not follow from tree convergence. In the following corollary we explore additional conditions in the tree convergence setting of Theorem \ref{prop7}, under which we also obtain residual mass convergence. Roughly speaking, condition (ii) below expresses the following intuition: we need label 1 in the asymptotically largest block most of the time, and on the corresponding set $\{|\Gamma|^\downarrow_1=|\Gamma_1|\}$ of infinite $\kappa$-measure, $|\Gamma^{[n]}|^\downarrow_1$ and $|\Gamma_1^{[n]}|$ approach their limit $|\Gamma|^\downarrow_1=|\Gamma_1|$ in a sufficiently regular way. The following statement includes Corollary \ref{cor p-con}.
\begin{cor} In the setting of Theorem \ref{prop7} {\rm(ii)}, the block $\Gamma_1$ containing {\rm 1} of $\kappa$-a.e.\ $\Gamma\in \cP$ has an asymptotic frequency in $(0,1)$. 
With $\Lambda$ as in Theorem \ref{propmc},
the following are equivalent:
  \begin{enumerate}\item[\rm(i)] $\displaystyle\int_{(0,\infty)}(1-e^{-x})\Lambda(dx)<\infty$ and $\displaystyle\int_{\cP}\left(|\Gamma_1^{[n]}| -|\Gamma_1|\right)\kappa(d\Gamma)\rightarrow 0$;
    \item[\rm(ii)] $\displaystyle\kappa(|\Gamma_1|\neq|\Gamma|^\downarrow_1)<\infty$ and $\displaystyle\int_{\{|\Gamma|^\downarrow_1=|\Gamma_1|\}}\left(|\Gamma^{[n]}|^\downarrow_1 -|\Gamma_1^{[n]}| \right)\kappa(d\Gamma)\rightarrow 0$.
  \end{enumerate}
  If condition {\rm(ii)} holds, then $X^{(n)}_{\lfloor\lambda_n t\rfloor}/n\!\rightarrow\! X_t$ in distribution, as $n\!\rightarrow\!\infty$, in the Skorohod sense as functions of $t\ge 0$, and this convergence holds jointly with the convergence of $\lambda_n^{-1} A_n$ to the absorption time of $X$ at $0$, where our notation is as in Theorem \ref{propmc}. 
\end{cor}
\begin{pf} Since $\Gamma$ has asymptotic ranked frequencies, $\Gamma_1$ has an asymptotic frequency by Lemma \ref{lm7}, for $\kappa$-a.e.\ $\Gamma\in\cP$. ``(i)$\Rightarrow$(ii)'' follows straight from the proof of Theorem \ref{propmc}, since (i) puts us into that setting; and also, the convergence of $X^{(n)}_{\lfloor\lambda_n t\rfloor}/n$ holds under (i). It remains to prove
``(ii)$\Rightarrow$(i)'', so we note that under (ii),
  \begin{equation}\label{upperbound}\int_{(0,\infty)}(1-e^{-x})\Lambda(dx)=\int_\cP\left(1-\left|\Gamma_1\right|\right)\kappa(d\Gamma)\le\kappa(\left|\Gamma_1\right|\neq\left|\Gamma\right|^\downarrow_1)+\int_\cP\left(1-\left|\Gamma\right|^\downarrow_1\right)\kappa(d\Gamma)<\infty
  \end{equation}
  and $\int_\cP\left(\left|\Gamma_1^{[n]}\right|-\left|\Gamma_1\right|\right)\kappa(d\Gamma)\rightarrow 0$ follows by the same argument as convergence in (\ref{split}), with roles of $(|\Gamma_1|,|\Gamma_1^{[n]}|)$ and $(|\Gamma|^\downarrow_1,|\Gamma^{[n]}|^\downarrow_1)$ interchanged, using Theorem \ref{prop7}(iii) for the first term and using the second condition under (ii) here for the third term of the modification of (\ref{split}). 
\end{pf}

\section{Further problems and related work}\label{further}

Due to the coupling of $(\cT_n,n\ge 1)$ in a regenerative tree growth process, the convergence in distribution in Theorems \ref{prop7} and \ref{propmc} should be strengthened to a convergence in probability or even to almost sure convergence in all cases discussed here. We have proved tree convergence in probability in the exchangeable case \cite{HMPW}, and in the restricted exchangeable case \cite{CW} provided that $\nu_j=\nu_m$, $j\ge m$, but the general case including the alpha-theta model remains open. 

In the alpha-theta model \cite{PW09} and the (restricted) exchangeable \cite{HMPW,CW} cases, we have established a two-stage almost sure convergence to a self-similar tree $\cT$ by passing via reduced subtrees of $\cT_n$ and of $\cT$ spanned by the first $k$ labelled leaves and letting first $n\rightarrow\infty$ and then $k\rightarrow\infty$. More specifically, we have embedded $(\cT_n,n\ge 1)$ in $\cT$ as discrete trees with edge lengths.

The basic embedding problem is to find a random leaf in a self-similar tree $(\cT,\mu)$ that induces a given decreasing self-similar Markov process as residual mass process, i.e. as the process that is parametrised by distance from the root on the path to the random leaf and that records for each point on the path the $\mu$-mass in the subtree above the point. Another interesting structure is the joint distribution of two residual mass processes (see \cite{PW09,PW13}).
When embedded in the same tree, they coincide up to a branch point and then evolve independently. In \cite{PW13}, we use the terms \em fragmenter \em for exponential subordinators $(e^{-\xi_s},s\ge 0)$, which are time-changed in Lamperti's representation 
(\ref{pssmp}), and \em bifurcator \em for pairs of fragmenters that coincide up to an exponential time and then evolve independently. In \cite{PW13} we investigate the fact that not all fragmenters appear as residual mass processes of typical (uniformly random) leaves. We introduce the notion of Markovian embedding in an exchangeable fragmentation process and show that for every (pure-jump) fragmenter $X$ there is a unique exchangeable dislocation measure $\kappa$ such that $X$ has a Markovian embedding into an associated exchangeable fragmentation process.  

In \cite{PW09,PW13}, we study an autonomous description of the evolution of reduced subtrees, viewed as weighted trees equipped with an (atomic) measure on the branches. We refer to a single branch with an atomic measure as a \em string of beads\em, see also \cite{Pal} for related structures. We refer to the evolution of reduced subtrees as \em bead splitting\em. In \cite{PW13}, we study certain binary bead splitting processes that evolve by size-biased branching, i.e.\ where an atom (a bead) is selected at random according to the measure on the branches and replaced by a (rescaled independent) copy of a given string of beads. We study the convergence of bead-splitting processes to self-similar CRTs.

\appendix
\section{Proof of Proposition \ref{regfromg}}\label{ap_regfromg}

First consider a regenerative tree growth rule, i.e.\ a sequence of transition 
  probability matrices $g_n$ from $\cP_n\setminus\{1_{[n]}\}$ to $\{0,\ldots,n+1\}$ with $g_n(\pi,0)$
  independent of $\pi$ and $g_n(\pi,i)=0$ if $\pi$ has strictly fewer than $i-1$ blocks. 
  For $n=1$ and $n=2$ the regenerative property is trivial. Consider the
  induction hypothesis that the growth rule gives rise to distributions $Q_m$ on $\bT_m$, 
  $m\le n$, and hence to $Q_B$ on $\bT_B$ after relabelling via the increasing bijection $[m]\rightarrow B$,
  for all $B\subset\bN$ with $\#B=m$, such that conditionally given a first split 
  $\Pi_n=(B_1,\ldots,B_k)$, the subtrees above the first split are independent, and the $i$th subtree 
  $\cT_{n,B_i}$ has conditional distribution 
  $Q_{B_i}$, $1\le i\le k$. For the induction step, note that conditionally given 
  $\Pi_n=(B_1,\ldots,B_k)$, the tree growth step from $n$ to $n+1$ specifies $Q_{n+1}$ on each of 
  the events $G_{n,i}$:
  \begin{itemize}\item $G_{n,0}$: here, $\Pi_{n+1}=([n],\{n+1\})$ is not related to $\Pi_n$; we will get back to this; 
    \item $G_{n,i}$, $1\le i\le k$: here, $\Pi_{n+1}=(B_1,\ldots,B_{i-1},B_i\cup\{n+1\},B_{i+1},\ldots,B_k)$; regenerative growth in the $i$th
      subtree preserves the conditional independence of subtrees, 
      and the induction hypothesis also yields that $\cT_{n+1,B_j}=\cT_{n,B_j}$ has
      conditional distribution $Q_{B_j}$ for $j\neq i$, while $\cT_{n+1,B_i\cup\{n+1\}}$ has conditional distribution
      $Q_{B_i\cup\{n+1\}}$ obtained from $Q_{B_i}$ via the growth rule applied to $B_i$ with
      $\#B_i\le n-1$. 
    \item $G_{n,k+1}$: here, $\Pi_{n+1}=(B_1,\ldots,B_k,\{n+1\})$, and conditional independence of
      subtrees as well as conditional distributions follow from the induction hypothesis, with the 
      addition of $\cT_{n+1,\{n+1\}}$ with (degenerate) conditional distribution $Q_{\{n+1\}}$.
  \end{itemize} 
  Conditionally given $\Pi_n=(B_1,\ldots,B_k)$, the events $G_{n,i}$ for $i\in\{k+2,\ldots,n+1\}$ 
  have probability zero, since $g_n(\pi,i)=0$ if $\pi$ has strictly fewer than $i-1$ blocks. Hence,
  $Q_{n+1}$ is fully specified and satisfies the regenerative property for each of the $i\ge 1$; 
  for the remaining $i=0$ case, we cannot work conditionally given $\Pi_n=(B_1,\ldots,B_k)$, because the regenerative
  property here is a statement conditionally given $\Pi_{n+1}=([n],\{n+1\})$, and indeed, since
  $g_n(\pi,0)$ does not depend on $\pi\in\cP_n\setminus\{1_{[n]}\}$, the subtree $\cT_{n+1,[n]}=\cT_n$ has
  conditional distribution $Q_{[n]}=Q_n$, while $\cT_{n+1,\{n+1\}}$ has conditional distribution
  $Q_{\{n+1\}}$. The induction proceeds.  

To prove the other direction, let $(\cT_n,n\ge 1)$ be a regenerative tree growth process. Consistency implies that the splitting rules $p_n(\pi)=\bP(\Pi_n=\pi)$ satisfy, for $\pi=(B_1,\ldots,B_k)\in\cP_n\setminus\{1_{[n]}\}$,
$$p_n(B_1,\ldots,B_k)=p_{n+1}([n],\{n+1\})p_n(B_1,\ldots,B_k)+\sum_{i=1}^{k+1}p_{n+1}(B_1,\ldots,B_i\cup\{n+1\},\ldots,B_k).
$$
  For $g_n(0)$ and $g_n(\pi,i)$ defined from $(p_n,n\ge 2)$ via \eqref{pg}, this implies that $g_n(\pi,i)\in[0,1]$ 
  and
  $$\sum_{i=0}^{k+1}g_n(\pi,i)=p_{n+1}([n],\{n+1\})+\sum_{i=1}^{k+1}\frac{p_{n+1}(B_1,\ldots,B_i\cup\{n+1\},\ldots,B_k)}{p_n(B_1,\ldots,B_k)}=1.
$$
  Also, \eqref{treeprob} 
  holds and determines $\bP(\cT_n=\ft\,|\,\cT_{n-1})$ as required, since $\cT_n$ determines $\cT_{n-1}$.\hfill $\square$

\section{Proof of Lemma \ref{lm7}}\label{ap_lemma7}
We consider the set\vspace{-0.1cm}
\[c_0 = \left\{ (s_1,s_2,\dots) \in [0,1]^\bN\colon \lim_{i\to\infty} s_i =0 \right\},\]
which is equipped with the uniform norm $||\cdot||_\infty$.  This set is clearly closed when considered as a subset of $\ell^\infty$ and thus is a complete metric space.  Let $F\colon c_0\to c_0$ be the map defined by $F(\mathbf{s}) = \mathbf{s}^\downarrow$, that is, $F$ is the map that takes a sequence to its non-increasing rearrangement.  Our first step is to prove that $F$ is continuous since this immediately implies that  if $|\Gamma^{[n]}|$ converges uniformly, say to $(y_i)_{i\geq 1}$, then $|\Gamma^{[n]}|^{\downarrow}$ converges to the non-increasing rearrangement of $(y_i)_{i\geq 1}$.

Fix $\epsilon>0$ and $\mathbf{s}\in c_0$.  Without loss of generality, we may assume that $\epsilon<\sup_i s_i$ and $\epsilon\notin\{s_i,i\ge 1\}$.  Let $B_{\bs}=\{i\ge 1\colon s_i> \epsilon\}$.  The fact that $\bs\in c_0$ implies that $ \#B_{\bs}<\infty$.  Observe that $F(\bs)$ is equal to the sequence obtained by concatenating the non-increasing rearrangement of $(s_i\colon i\in B_{\bs})$ with the non-increasing rearrangement of $(s_i\colon i\notin B_{\bs})$. 

Suppose that $\mathbf{s}_n\rightarrow \mathbf{s}$.  For sufficiently large $n$ we have $B_{\bs_n}=B_{\bs}$.  Since ranking is continuous on $\bR^{\#B_\fs}$, it follows that\vspace{-0.1cm} 
\[\lim_{n \to \infty} (s^\downarrow_{n,1},\dots, s^\downarrow_{n,\#B_{\bs}})  = (s^\downarrow_1,\dots,s^\downarrow_{\#B_{\bs}})\] 
and also that\vspace{-0.1cm}
\[\sup_{i>\# B_{\bs}} s^\downarrow_i+\limsup_{n\to\infty}\sup_{i>\#B_{\bs}} s^{\downarrow}_{n,i} \leq 2\epsilon.\] 
As a result we have
\[\limsup_{n\to\infty} ||F(\bs)-F(\bs_n)||_\infty \leq 2\epsilon,\]
and the continuity of $F$ follows.\pagebreak[1] 

We now prove the opposite direction.  To that end, assume that  $|\Gamma^{[n]}|^{\downarrow} \to (x_i)_{i\geq 1}$ pointwise.  We will prove that $|\Gamma^{[n]}|$ converges in $c_0$ and the proof of the previous part then identifies the limit.  Since $|\Gamma^{[n]}|^{\downarrow}$ is non-increasing for each $n$ with sums uniformly bounded by $1$, this implies that $|\Gamma^{[n]}|^{\downarrow} \to (x_i)_{i\geq 1}$ uniformly.  If $x_1=0$ we are done, so we assume that $x_1>0$.  Let $\epsilon>0$ be given, and without loss of generality suppose that $\epsilon <x_1$.  By Fatou's lemma we have $\sum_{i\geq 1} x_i \leq 1 <\infty$ and, consequently, we can choose $K$ so that $\sum_{i\geq K+1} x_i < \epsilon$.  Let 
\[ \epsilon_1 = \epsilon \wedge \min\left\{ \frac{|x_i-x_j|}{3}\colon 1\leq i,j\leq K+1 \textrm{ and } x_i\neq x_j\right\}.\]
Since $|\Gamma^{[n]}|^{\downarrow} \to (x_i)_{i\geq 1}$ uniformly, we can choose $N>1/\epsilon_1$ such that $\sup_i \left||\Gamma^{[n]}|_i^{\downarrow} - x_i\right|< \epsilon_1$ for all $n\geq N$.  For each $n\geq N$ let $\sigma_n\colon\bN\to \bN$ be a bijection such that $(|\Gamma^{[n]}|_{\sigma_n(i)})_{i\geq 1} = |\Gamma^{[n]}|^{\downarrow}$.  Note that we have used the fact that $|\Gamma^{[n]}|$ has only finitely many non-zero entries to obtain this bijection.  Since $N>1/\epsilon_1$, for all $i\ge 1$ and $n\geq N$ we have 
\[ \left| |\Gamma^{[n]}|_i - |\Gamma^{[n+1]}|_i\right| \leq 1/(n+1) < \epsilon_1.\]
It follows that for $n\geq N$\vspace{-0.1cm} 
\[ \sup_i\left| |\Gamma^{[n+1]}|_{\sigma_n(i)} - x_i\right| \leq 2\epsilon_1.\]
By our choice of $\epsilon_1$, for any $1\leq j \leq K$ and any $i\ge 1$ such that $x_j\neq x_i$ we have
\[ \left||\Gamma^{[n+1]}|_{\sigma_n(i)} - x_j\right| \geq  |x_i-x_j| - \left||\Gamma^{[n+1]}|_{\sigma_n(i)} - x_i\right| \geq  \epsilon_1.\]\pagebreak

\noindent However, since $\sup_i \left||\Gamma^{[n]}|_i^{\downarrow} - x_i\right|< \epsilon_1$ for all $n\geq N$, this implies that
\[  \sup_{1\leq i\leq K}\left| |\Gamma^{[n+1]}|_{\sigma_n(i)} - x_i\right| < \epsilon_1 \quad \textrm{and}\quad \sup_{i\geq K+1} |\Gamma^{[n+1]}|_{\sigma_n(i)} < x_{K+1} +\epsilon_1. \]
Inductively, we conclude that for all $n\geq N$ and $k\geq 0$ 
$$ \sup_{1\leq i\leq K}\left| |\Gamma^{[n+k]}|_{\sigma_n(i)} - x_i\right| < \epsilon_1 \quad \textrm{and}\quad \sup_{i\geq K+1} |\Gamma^{[n+k]}|_{\sigma_n(i)} < x_{K+1} +\epsilon_1.$$
Combining these, we see that for $n\geq N$
$$ \sup_{i\geq 1} \left| |\Gamma^{[n]}|_{i} - x_{\sigma_N^{-1}(i)}\right| = \sup_{i\geq 1} \left| |\Gamma^{[n]}|_{\sigma_N(i)} - x_i\right|  < 2\epsilon.$$
We are not quite done since $\sigma_N$ depends on $\epsilon_1$.  Note, however, that the above inequality implies for $n, m \geq N$
\[ \sup_{i\geq 1} \left| |\Gamma^{[n]}|_{i} -|\Gamma^{[m]}|_{i} \right| < 4\epsilon.\] 
This shows that $(|\Gamma^{[n]}|)_{n\geq 1}$ is a Cauchy sequence in the complete metric space $c_0$ and, therefore, converges uniformly.\hfill $\square$

\bibliographystyle{abbrv}
\bibliography{genreg}

\begin{thebibliography}{10}

\bibitem{Ald-91}
D.~Aldous.
\newblock The continuum random tree. {I}.
\newblock {\em Ann. Probab.}, 19(1):1--28, 1991.

\bibitem{Aldo93}
D.~Aldous.
\newblock The continuum random tree. {III}.
\newblock {\em Ann. Probab.}, 21(1):248--289, 1993.

\bibitem{Ald-93}
D.~Aldous.
\newblock Probability distributions on cladograms.
\newblock In {\em Random discrete structures (Minneapolis, MN, 1993)},
  volume~76 of {\em IMA Vol. Math. Appl.}, pages 1--18. Springer, New York,
  1996.

\bibitem{Ber-book}
J.~Bertoin.
\newblock {\em Random fragmentation and coagulation processes}, volume 102 of
  {\em Cambridge Studies in Advanced Mathematics}.
\newblock Cambridge University Press, Cambridge, 2006.

\bibitem{CFW}
B.~Chen, D.~Ford, and M.~Winkel.
\newblock A new family of {M}arkov branching trees: the alpha-gamma model.
\newblock {\em Electron. J. Probab.}, 14:no. 15, 400--430 (electronic), 2009.

\bibitem{CW}
B.~Chen and M.~Winkel.
\newblock Restricted exchangeable partitions and embedding of associated
  hierarchies in continuum random trees.
\newblock {\em Ann. Inst. Henri Poincare Probab. Stat.}, 49(3):839--872, 2013.

\bibitem{Du-LG}
T.~Duquesne and J.-F. Le~Gall.
\newblock Random trees, {L}\'evy processes and spatial branching processes.
\newblock {\em Ast\'erisque}, (281):vi+147, 2002.

\bibitem{EPW}
S.~N. Evans, J.~Pitman, and A.~Winter.
\newblock Rayleigh processes, real trees, and root growth with re-grafting.
\newblock {\em Probab. Theory Related Fields}, 134(1):81--126, 2006.

\bibitem{EW06}
S.~N. Evans and A.~Winter.
\newblock Subtree prune and regraft: a reversible real tree-valued {M}arkov
  process.
\newblock {\em Ann. Probab.}, 34(3):918--961, 2006.

\bibitem{For-05}
D.~J. Ford.
\newblock Probabilities on cladograms: introduction to the alpha model.
\newblock {\em arXiv:math.PR/0511246}, 2005.

\bibitem{gp}
A.~Gnedin and J.~Pitman.
\newblock Regenerative composition structures.
\newblock {\em Ann. Probab.}, 33(2):445--479, 2005.

\bibitem{GPW09}
A.~Greven, P.~Pfaffelhuber, and A.~Winter.
\newblock Convergence in distribution of random metric measure spaces
  ({$\Lambda$}-coalescent measure trees).
\newblock {\em Probab. Theory Related Fields}, 145(1-2):285--322, 2009.

\bibitem{G99}
M.~Gromov.
\newblock {\em Metric structures for {R}iemannian and non-{R}iemannian spaces},
  volume 152 of {\em Progress in Mathematics}.
\newblock Birkh\"auser Boston Inc., Boston, MA, 1999.
\newblock Based on the 1981 French original, with appendices by M. Katz, P.
  Pansu and S. Semmes, Translated from the French by Sean Michael Bates.

\bibitem{HM04}
B.~Haas and G.~Miermont.
\newblock The genealogy of self-similar fragmentations with negative index as a
  continuum random tree.
\newblock {\em Electron. J. Probab.}, 9:no. 4, 57--97 (electronic), 2004.

\bibitem{HM09}
B.~Haas and G.~Miermont.
\newblock {Self-similar scaling limits of non-increasing Markov chains}.
\newblock {\em Bernoulli}, 17(4):1217--1247, 2011.

\bibitem{HM10}
B.~Haas and G.~Miermont.
\newblock {Scaling limits of Markov branching trees, with applications to
  Galton-Watson and random unordered trees}.
\newblock {\em Ann. Probab.}, 40(6):2589--2666, 2012.

\bibitem{HMPW}
B.~Haas, G.~Miermont, J.~Pitman, and M.~Winkel.
\newblock Continuum tree asymptotics of discrete fragmentations and
  applications to phylogenetic models.
\newblock {\em Ann. Probab.}, 36(5):1790--1837, 2008.

\bibitem{HPW}
B.~Haas, J.~Pitman, and M.~Winkel.
\newblock Spinal partitions and invariance under re-rooting of continuum random
  trees.
\newblock {\em Ann. Probab.}, 37(4):1381--1411, 2009.

\bibitem{HaP-11}
C.~Haulk and J.~Pitman.
\newblock A representation of exchangeable hierarchies by sampling from real
  trees.
\newblock {\em arXiv:1101.5619 [math.PR]}, 2011.

\bibitem{Lam-72}
J.~Lamperti.
\newblock Semi-stable {M}arkov processes. {I}.
\newblock {\em Z. Wahrscheinlichkeitstheorie und Verw. Gebiete}, 22:205--225,
  1972.

\bibitem{Mar-08}
P.~Marchal.
\newblock A note on the fragmentation of a stable tree.
\newblock In {\em Fifth Colloquium on Mathematics and Computer Science},
  volume~AI, pages 489--500. Discrete Mathematics and Theoretical Computer
  Science, 2008.

\bibitem{MPW}
P.~McCullagh, J.~Pitman, and M.~Winkel.
\newblock Gibbs fragmentation trees.
\newblock {\em Bernoulli}, 14(4):988--1002, 2008.

\bibitem{Mie-03}
G.~Miermont.
\newblock Self-similar fragmentations derived from the stable tree. {I}.
  {S}plitting at heights.
\newblock {\em Probab. Theory Related Fields}, 127(3):423--454, 2003.

\bibitem{M09}
G.~Miermont.
\newblock Tessellations of random maps of arbitrary genus.
\newblock {\em Ann. Sci. \'Ec. Norm. Sup\'er. (4)}, 42(5):725--781, 2009.

\bibitem{Pal}
S.~Pal.
\newblock {On the Aldous diffusion on Continuum Trees. I}.
\newblock {\em arXiv:1104.4186 [math.PR]}, 2011.

\bibitem{Pit-95}
J.~Pitman.
\newblock Exchangeable and partially exchangeable random partitions.
\newblock {\em Probab. Theory Related Fields}, 102(2):145--158, 1995.

\bibitem{PR}
J.~Pitman and D.~Rizzolo.
\newblock {Schr\"oder's problems and scaling limits of random trees}.
\newblock {\em Trans. Amer. Math. Soc., to appear, preprint available at
  arXiv:1107.1760 [math.PR]}, 2013.

\bibitem{PW09}
J.~Pitman and M.~Winkel.
\newblock Regenerative tree growth: binary self-similar continuum random trees
  and {P}oisson-{D}irichlet compositions.
\newblock {\em Ann. Probab.}, 37(5):1999--2042, 2009.

\bibitem{PW13}
J.~Pitman and M.~Winkel.
\newblock {Regenerative tree growth: Markovian embedding of fragmenters,
  bifurcators and bead splitting processes}.
\newblock {\em arXiv:1304.0802 [math.PR]}, 2013.

\end{thebibliography}

\newpage

\end{document}